\documentclass[a4paper]{article}
 
 \usepackage{a4wide}
\usepackage{latexsym}   

\usepackage{amsmath,amssymb,amsthm}
\usepackage{dsfont}
\usepackage{graphicx}
\usepackage[british]{babel}
\usepackage[utf8]{inputenc}
\usepackage{xcolor}
\usepackage[colorlinks=true,linkcolor=blue,citecolor=red]{hyperref}
\usepackage{authblk}
\usepackage{blindtext}
\usepackage{epstopdf}
\usepackage{caption}
\usepackage{subcaption}

\newtheorem{remark}{Remark}

\newcommand{\be}{\begin{equation}}
\newcommand{\ee}{\end{equation}}
\newcommand{\benn}{\begin{equation*}}
\newcommand{\eenn}{\end{equation*}}

\usepackage[margin=1cm]{caption}

\newcommand{\z}{\mathbf{z}}
\newcommand{\f}{\mathbf{f}}
\theoremstyle{plain}

\newtheorem{theorem}{Theorem}[section]

\DeclareMathOperator{\sign}{sign}

\begin{document}
\title{Uncertainty quantification and control of kinetic models {of} tumour growth under clinical uncertainties}

\author[5]{A. Medaglia}
\author[1,5]{G. Colelli }
\author[1]{L. Farina}
\author[1]{A. Bacila}
\author[2]{P. Bini}
\author[2]{E. Marchioni}
\author[4]{S. Figini}
\author[1,3]{A. Pichiecchio}
\author[5]{M. Zanella\thanks{\tt mattia.zanella@unipv.it}}

\affil[1]{Neuroradiology Department, IRCCS Mondino Foundation, Pavia, Italy}
\affil[2]{Neurology Department, IRCCS Mondino Foundation, Pavia, Italy }
\affil[3]{Department of Brain and Behavioural Sciences, University of Pavia, Italy}
\affil[4]{Department of Political and Social Sciences, University of Pavia, Italy}
\affil[5]{Department of Mathematics "F. Casorati", University of Pavia, Italy.}

\date{}

\maketitle

\abstract{
In this work, we develop a kinetic model {of} tumour growth taking into account the effects of clinical uncertainties characterising the tumours' progression. The action of therapeutic protocols trying to steer the tumours' volume towards a target size is then investigated by means of suitable selective-type controls acting at the level of cellular dynamics. By means of classical tools of statistical mechanics for many-agent systems, we are able to prove that it is possible to dampen clinical uncertainties across the scales. To take into account the scarcity of clinical data and the possible source of error in the image segmentation of tumours' evolution, we estimated empirical distributions of relevant parameters that are considered to calibrate the resulting model obtained from real cases of primary glioblastoma. Suitable numerical methods for uncertainty quantification of the resulting kinetic equations are discussed and, in the last part of the paper, we compare the effectiveness of the introduced control approaches in reducing the variability in tumours' size due to the presence of uncertain quantities. 
}
\\[+.2cm]

\tableofcontents

\section{Introduction}

The processes of tumour formation are highly complex phenomena involving different stages starting from damages in the DNA molecules leading to harmful mutations in the cell's genome that are not repaired in absence of cellular apoptosis. This mechanism leads to an unregulated mitosis and then to the formation of tumours. These mutations may be triggered by many aspects, including both environmental and genetic factors, see \cite{Frank,GCI,Folk}. 

In the last decades, extensive research efforts have been devoted to the mathematical formalisation of tumour growth dynamics and to the formalisation of growth factors, see e.g. \cite{AG,Gerlee,HCMB,Leroy,NR,RBKW,RCM}. Among the easiest way to describe these biological phenomena can be found in the literature of population dynamics to describe the evolution in time of the volume of a tumour. This modelling approach is based on first order ODEs that are named according to the form of the right-hand side. Leading examples are Gompertz and von Bertalanffy models. More recently, West and collaborators proposed a variation to the von Bertalanffy model in \cite{West}. It is worth mentioning that there is no unanimous consensus on the most appropriate modelling setting and several proposals have been introduced to better describe these dynamics. In particular, in \cite{PTZ} the authors proposed a statistical approach based on kinetic theory to describe the growth of tumour cells in terms of the evolution of a distribution function. The temporal variation of such distribution is considered as a result of elementary transitions occurring at the cellular level which takes environmental cues and random fluctuations into account. The expected cellular variations are coherent with the mentioned ODE-based models in suitable limits. Furthermore, through the explicit computation of the equilibrium states of the resulting Fokker-Planck-type equation, we get additional information on the decay of the tails. In particular, it is shown that von Bertalanffy-type models lead to fat tailed distributions of the volumes of tumours, whereas Gompertz-type models are linked to slim tailed distributions. The mathematical understanding of the behaviour of the tails is essential to quantify the probability of having tumours growing to sizes that are harmful to the human body. Existing kinetic models for statistical growth dynamics are linked to cell mutations \cite{KP,Tosc2}. {In particular, in recent years a huge literature on mathematical modelling for glioma growth have been developed, see \cite{Conte,EHKS,Painter-Hillen} and the references therein.  }

Even if the mathematical simplicity of ODE-based modelling allows to handle more efficiently parameter estimation tasks, see e.g. \cite{Ben,MVPF,Norton,WN,Whel,WK}, the models based on partial differential equations are capable to describe the phenomenon under study in a statistical way \cite{LP,PTZ} or highlighting the mechanical properties of the tissues, see e.g. \cite{Ag,GP}. Furthermore, the lack of accurate clinical data introduces many sources of uncertainties stemming out at various levels of observation when facing the progression of human cancer. To mention a few, the main limitation consists in a limited set of observations due to clinical constraints. The second one may arise from the manual corrections of 3D semi-automatic tumours segmentation. The third comes from the fact that the evolution may differ strongly from one individual to another, since in each host the response of the body is influenced by many factors, like the clinical history of a patient. For these reasons, to produce effective predictions and to better understand the physical phenomena under study, we incorporate ineradicable uncertainties in the dynamics from the beginning of the modelling. The introduction of uncertain quantities points in the direction of a more realistic description of the underlying processes and helps us to compute possible deviations from the prescribed deterministic behaviour.  

Once established the emerging distribution of the kinetic model in presence of uncertain quantities we further propose a robust approach to steer the system towards a prescribed target to mimic implementable therapeutic protocols. The control is here conceived as an additional external dynamics depending on the state of the system. The proposed control setting has roots in Boltzmann-type controls proposed in \cite{AFK,AP,Albi0,Albi1} where an optimal control problem is solved at the microscopic level and then studied at the mesoscopic scale through classical methods of kinetic theory  \cite{Cer,PT}. This approach has connections with classical approaches for the control of mean-field equations, see \cite{BFY}. Aside from the mentioned methods, the control of emergent behaviour has been studied also on the level of the microscopic agents, see e.g.  \cite{BBCK}, as well as fluid–dynamic equations. The contributions have to be further distinguished depending on the type of applied control. Without intending to review all literature we give some references on certain classes of control, e.g. sparse control \cite{FPR}, Nash equilibrium control \cite{DHL}, control using linearised dynamics and Riccati equations \cite{HSP}. 

In the proposed setting, we discuss analytical properties of the asymptotic regime highlighting that a damping of structural uncertainties of the system is achieved at the macroscopic/observable level. Furthermore, the proposed approach is genuinely multiscale since it makes it possible to bridge actions on the individual cellular-based dynamics to observable patterns in the cohort of patients. In a different context, the asymptotic properties of such controls have been investigated in \cite{TZ}.

From the mathematical viewpoint, the introduction of such clinical uncertainties translates in an increased dimensionality of the resulting kinetic problem whose equilibrium depends on all the uncertainties introduced at the cellular level. The construction of numerical schemes for the resulting equations needs to guarantee spectral convergence on the random field under suitable regularity assumptions together with the preservation of the main physical properties of the model, see e.g. \cite{CZ,CPZ,ParZan,Xiu}. 
 
In more details, the paper is organised as follows: in Section \ref{sect:kin} we introduce the kinetic model of interest and we discuss the role of the uncertain parameters present at the level of the transition function. Hence, we briefly derive in the quasi-invariant limit reduced order models of Fokker-Planck-type from which large time distributions are explicitly computable. In Section \ref{sect:ther} we introduce a hierarchical control protocol with the aim to steer the tumour's size towards a prescribed size through the minimisation of two possible cost functionals. The emerging macroscopic properties of the introduced approach is then discussed together with their interplay with the model uncertainties. In Section \ref{sec:QoI} we face the calibration of the model with real clinical data provided and in Section \ref{sect:numerics} we introduce suitable numerical strategies to deal with uncertainty quantification of Boltzmann-type and Fokker-Planck-type equations.

\section{Kinetic modelling of tumour growth dynamics with clinical uncertainties}
\label{sect:kin}

Let us characterise the microscopic state of an evolving tumour by means of a variable $x \in \mathbb R_+$ representing the volume of the tumour. Furthermore, we collect all the {sources} of uncertainties of the dynamics in the random vector $\z = (z_1,\dots,z_d)\in \mathbb R^{d}$ whose distribution is $\rho(\z)$, i.e.
\[
\mathbb P[\z \in A] = \int_A \rho(\z)dz,
\] 
for any $A \subseteq \mathbb R^d$.  In details, for any fixed $\z$, if $X(\z,t)$ is a random variable expressing the volume of the tumour, the probability density associated to $X(\z,t)$ is $f(\z,x,t)$ and $f(\z,x,t)dx$ is the fraction of tumours that, at time $t\ge 0$, are characterised by a volume between $x$ and $x+dx$. It is worth to mention that the knowledge of the evolution of $f(\z,x,t)$ allows to compute the evolution of macroscopic quantities that are given by 
\[
\int_{\mathbb R_+}\varphi(x)f(\z,x,t) dx, 
\]
where $\varphi$ is a test function. We can observe that the macroscopic quantity of interest still depends on the introduced uncertainties.

In details, for a given volume $x \in \mathbb R_+$ of cancer cells, we characterise an elementary variation $x \rightarrow x'$ as follows
\be
\label{eq:trans}
x' = x + \Phi^\epsilon_\delta(x/x_L,\z)x + x \eta_\epsilon, \qquad \epsilon\ll 1.
\ee
where  $\Phi^\epsilon_\delta$ is a transition function taking into account variations due to environmental cues and which depends on the tumour size $x$ and on additional clinical uncertainties expressed by the random vector $\z \in  \mathbb R^{d}$ characterising the lack of knowledge of parameters. The quantity $x_L = x_L(\z)>0$ is a characteristic patient-based tumour size, e.g. the carrying capacity. Furthermore, in \eqref{eq:trans} the random variable $\eta_\epsilon$ takes into account unpredictable changes in the transition dynamics and such that $\langle \eta_\epsilon \rangle = 0$ and $\langle \eta_{\epsilon}^2 \rangle = \epsilon\sigma^2$, where $\langle \cdot \rangle$ denotes the expectation with respect to the distribution of $\eta_\epsilon$. Therefore, in a single transition the tumour's size can be modified by two multiplicative mechanisms parametrised by the positive constant $\epsilon\ll 1$ and {by} the uncertain parameter $\delta = \delta(\z) \in [-1,1]$ influencing the considered type of growth.

\subsection{Transition functions}\label{sect:transition}

The transitions characterising the proposed elementary growths should be considered arbitrary small. For this reason, coherently with \cite{PTZ},  we require that $\Phi^\epsilon_\delta$ is of order $\epsilon$ and that 
\[
\lim_{\epsilon \rightarrow 0^+}\dfrac{\Phi^\epsilon_\delta(x/x_L,\z)}{\epsilon} = \Phi_\delta(x/x_L,\z). 
\]
Having in mind this requirement we now characterise a general transition function that is coherent with known microscopic models for tumour growth. We consider 
\be
\label{eq:Phi}
\Phi^\epsilon_\delta(y,\z) = \mu \dfrac{1-e^{\epsilon(y^{\delta}-1)/\delta}}{(1+\lambda)e^{\epsilon(y^{\delta}-1)/\delta}+1-\lambda},\qquad y = \frac{x}{x_L}
\ee
where we introduced the uncertain parameters $\mu = \mu(\z) \in (0,1)$ and $\lambda = \lambda(\z) \in [0,1)$ characterising birth and death rates of tumour cells in a single transition since, independently on $\epsilon\ll 1$, we have
\[
- \dfrac{\mu}{1+\lambda} \le \Phi^\epsilon(x/x_L,\z) \le \dfrac{\mu}{1-\lambda}. 
\]
{In absence of fluctuations, we have $x^\prime>x$ when $x < x_L$ for all values of the parameter $\delta$. In terms of $\delta$, the transition function do not behave in the same way in the region $x < x_L$. As highlighted in \cite{DT20,PTZ}, the transition function $\Phi_\delta^\epsilon$ with $\delta>0$ is increasing convex for all $x \le x_L$, whereas it is concave in an interval $[0,\bar x]$, $\bar x <x_L$ and then convex for $\delta<0$}
A compatibility condition for the elementary variations \eqref{eq:trans} with the transition functions \eqref{eq:Phi} is that $x^\prime$ remains in $\mathbb R_+$. This can be guaranteed by imposing the following sufficient condition on the fluctuation $\eta_\epsilon$. In particular, by considering for any $\z \in \mathbb R^d$ a random variable such that
\[
\eta_\epsilon \ge -1 + \max_{\z \in \mathbb R^d} \dfrac{\mu}{1+\lambda},
\]
the post-transition size  $x^\prime$ is positive. 

It is worth to remark that in the limit $\epsilon \rightarrow 0^+$ we have
\[
\Phi^\epsilon(x/x_L,\z) \approx \epsilon\mu \dfrac{(y^{\delta}-1)/\delta}{(1+\lambda)\epsilon (1-y^{\delta})/\delta+2},\qquad y = \frac{x}{x_L}
\] 
which implies 
\[
\lim_{\epsilon \rightarrow 0^+}\dfrac{\Phi^\epsilon(x/x_L,\z)}{\epsilon} = \dfrac{\mu}{2\delta}\left( 1-\left(\dfrac{x}{x_L} \right)^{\delta}\right).
\]
Therefore, the proposed transition function is coherent in the limit $\epsilon \rightarrow 0^+$ with existing models for the description of tumour growth. In particular,  if we consider the following first order Bernoulli differential equation
\be\label{eq:micro}
\dfrac{d}{dt}x(\z,t) = \dfrac{\mu}{2\delta} \left( 1- \left( \dfrac{x(\z,t)}{x_L} \right)^{\delta} \right) x(\z,t),
\ee
in the limit $\delta \rightarrow 0$ we recover Gompert growth dynamics since \eqref{eq:micro} corresponds to
\begin{equation}
\label{eq:gomp}
\dfrac{d}{dt}x(\z,t) = 
-\dfrac{\mu}{2} \log\left(\dfrac{x(\z,t)}{x_L} \right)x(\z,t),
\end{equation}
whereas for any $\delta<0$ we recover von Bertalanffy dynamics of the form 
\begin{equation}
\label{eq:vb}
\dfrac{d}{dt}x(\z,t) = px(\z,t)^{\delta+1} - q x(\z,t), 
\end{equation}
with $q = q(\z) = -\dfrac{\mu}{2\delta}$, $p = p(\z) = -\dfrac{\mu}{2\delta x_L^\delta}$. It can be easily observed that for any $\delta >0$ we recover logistic-type growth that are not of limited interest in the present context. 

\begin{remark}
The dynamics described by \eqref{eq:micro} are coherent with the expected transition scheme \eqref{eq:trans}. Indeed, if we introduce a forward time discretisation with time step $\Delta t>0$ from \eqref{eq:micro} we get
\[
x^{n+1}(\z) = x^n(\z) +  \Delta t \dfrac{\mu}{2\delta} \left( 1- \left( \dfrac{x^n(\z)}{x_L} \right)^{\delta} \right) x^n(\z),
\]
where $x^n(\z) = x(t^n,\z)$ and $t^n = n \Delta t$, $n \in \mathbb N$. Hence, by identifying $\Delta t = \epsilon$ we can recognise the transition scheme in \eqref{eq:trans}. 
\end{remark}

\subsection{Kinetic models and equilibria}
\label{subsec:equilibria}

Let $f = f(\z,x,t)$ be the distribution function of cells of size $x \in \mathbb R_+$ at time $t \ge0$  and depending on the epistemic uncertainties collected in $\z$. The evolution of $f$ is the given by the following kinetic equation
\be\label{eq:kin_strong}
\partial_t f(\z,x,t) =  Q_G(f)(\z,x,t), 
\ee
where the transition operator $Q_G(\cdot)(\z,x,t)$ is defined as follows
\be\label{eq:QG}
Q_G(f)(\z,x,t) = \int_{\mathbb R_+} \dfrac{1}{{}^\prime \mathcal J}f(\z,{}^\prime x,t)dx-f(\z,x,t), 
\ee
being ${}^\prime x$ the pre-transition state {and ${}^\prime\mathcal J$ is the absolute value of the Jacobian of the transformation from the pre-transition state ${}^\prime x$ to the state  $x$}. The kinetic equation \eqref{eq:kin_strong} can be fruitfully written in weak form to evaluate the evolution of observable quantities as follows
\be\label{eq:kinetic}
\dfrac{d}{dt} \int_{\mathbb R_+} \varphi(x)f(\z,x,t)dx =  \int_{\mathbb R_+}\langle\varphi(x^\prime)-\varphi(x)\rangle f(\z,x,t)dx, 
\ee
where $\varphi$ is a smooth function. Since the {analytical} computation of {the} equilibrium distribution of \eqref{eq:kin_strong} is very hard it is convenient to resort to a surrogate model  with which we can analytically obtain the large type distribution of the studied phenomenon. This approach is defined as quasi-invariant limit and it has roots in the grazing limit of kinetic theory. Several applications of this approach have been employed in recent years for the statistical description of collective phenomena, see \cite{PT,Tosc1} for an introduction. In the following we briefly recall the derivation of Fokker-Planck-type equations from \eqref{eq:kin_strong} thanks to a quasi-invariant limit technique. 

We may observe that for $\epsilon \ll 1$ the difference $x^\prime-x$ is small and we can perform the following Taylor expansion up to order three
\[
\varphi(x^\prime)-\varphi(x) = (x^\prime -x)\dfrac{d \varphi(x)}{dx} + \dfrac{1}{2} (x^\prime-x)^2  \dfrac{d^2\varphi(x)}{dx^2} + \dfrac{1}{6} (x^\prime-x)^3  \dfrac{d^3\varphi(\bar x)}{dx^3},
\] 
with $\bar x \in (\min\{x,x^\prime\},\max\{x,x^\prime\})$. Since $x^\prime-x = \Phi^\epsilon_\delta(x/x_L,\z) + x \eta_\epsilon$ we can plug this expression in \eqref{eq:kinetic} to obtain 
\be \label{eq:boltz}
\begin{split}
    \dfrac{d}{dt} \int_{\mathbb R_+}\varphi(x)f(\z,x,t)dx =& \int_{\mathbb R_+} \dfrac{\Phi^{\epsilon}_\delta(x/x_L,\z)}{\epsilon}x f(\z,x,t)\dfrac{d\varphi(x)}{dx}dx + \dfrac{\sigma^2}{2} \int_{\mathbb R_+}x^2f(\z,x,t)\dfrac{d^2\varphi(x)}{dx^2}dx \\
    & + R_{\varphi}(f)(\z,x,t),
\end{split}
\ee
where we have exploited the fact that $\langle \eta_{\epsilon} \rangle = 0$ and we have defined the rest $R_{\varphi}(f)(\z,x,t)$ as
\benn
\begin{split}
    R_{\varphi}(f)(\z,x,t) := &  \dfrac{1}{2}\int_{\mathbb R_+} \dfrac{(\Phi^{\epsilon}_\delta(x/x_L,\z))^2}{\epsilon}x^2 
    f(\z,x,t)\dfrac{d^2\varphi(x)}{dx^2}dx \\
    & + \dfrac{1}{6} \int_{\mathbb R_+} \dfrac{\langle \Phi^\epsilon_\delta(x/x_L,\z)x + x \eta_\epsilon \rangle^3}{\epsilon} f(\z,x,t)\dfrac{d^3\varphi(x)}{dx^3}dx.
\end{split}
\eenn
Assuming that the third order {moment} of $\eta_\epsilon$ is bounded, i.e. $\langle |\eta_\epsilon|^3 \rangle < +\infty$, thanks to the smoothness of $\varphi$ we have 
\[
|R_{\varphi}(f)(x,t,\z)| \lesssim \epsilon + \epsilon^2 + \epsilon + \sqrt{\epsilon},
\]
where we use the notation $a\lesssim b$ to mean that there exists a constant $K>0$ such that $a\leq Kb$. Hence, in the limit $\epsilon\to0^+$ we have $|R_{\varphi}(f)(\z,x,t)|\to 0$, for every $x\in\mathbb R_+$, $t>0$ and $\z \in \mathbb R^d$. As a consequence, if we introduce the new time scale $\tau = \epsilon t$, for $\epsilon\ll 1$, {a distribution function $f(\z,x,\tau) = f(\z,x,\tau/\epsilon)$,} and we observe that $\frac{d}{dt} = \epsilon \frac{d}{d\tau}$, the model defined in \eqref{eq:boltz} for $\epsilon \rightarrow 0^+$ converges to
\[
\dfrac{d}{d\tau} \int_{\mathbb R_+}\varphi(x)f(\z,x,\tau)dx = \int_{\mathbb R_+} \Phi_\delta(x/x_L,\z)x f(\z,x,\tau)\dfrac{d\varphi(x)}{dx}dx + \dfrac{\sigma^2}{2} \int_{\mathbb R_+}x^2f(\z,x,\tau)\dfrac{d^2\varphi(x)}{dx^2}dx,
\]
where
\be\label{eq:Phi_res}
\Phi_\delta(x/x_L,\z) = \dfrac{\mu}{2\delta}\left(1- \left(\dfrac{x}{x_L} \right)^{\delta} \right).
\ee
Hence, integrating back by parts we obtain the following Fokker-Planck-type equation with uncertainties
\be
\label{eq:FP}
\partial_\tau f(\z,x,\tau) = \partial_x \left[ -\Phi_\delta(x/x_L,\z) x f(\z,x,\tau) + \dfrac{\sigma^2}{2}\partial_x (x^2 f(\z,x,\tau) )\right].
\ee
provided that for all $\tau\ge 0$ the density $f(\z,x,\tau)$ satisfies the following no-flux boundary condition
\begin{equation}
\label{eq:BC}
-\Phi_\delta(x/x_L,\z)x f(\z,x,\tau) + \dfrac{\sigma^2}{2} \partial_{x}(x^2 f(\z,x,\tau)) \Bigg|_{x=0}=0.
\end{equation}

Thanks to the obtained surrogate model we can study the large time behaviour of the system. {In particular, the model \eqref{eq:FP} with no-flux boundary condition \eqref{eq:BC}  admits a unique equilibrium distribution $f^\infty(\z,x)$ that is solution to
\[
 -\Phi_\delta(x/x_L,\z) x f^\infty(\z,x) + \dfrac{\sigma^2}{2}\partial_x (x^2 f^\infty(\z,x) ) = 0, 
\]
see \cite{Risken}.}
In view of \eqref{eq:Phi_res} we have
\be \label{eq:equilibria}
f^\infty(\z,x) =C_{\mu,\sigma^2,x_L}(\z) \left( \dfrac{x}{x_L} \right)^{\frac{\mu}{\sigma^2 \delta}-2} \exp\left\{ -\dfrac{\mu}{\sigma^2 \delta} \left( \left( \dfrac{x}{x_L} \right)^{\delta}-1 \right) \right\},
\ee
with $C_{\mu,\sigma^2,x_L}>0$ a normalisation constant. 

In particular, we highlight that the two reference microscopic models we consider, corresponding to the choices $\delta<0$ and $\delta \rightarrow 0$, generate slight different equilibria. In particular, the Gompertzian growths, obtained in the limit  $\delta \rightarrow 0$, generates  at the equilibrium the lognormal distribution
\[
f^\infty(\z,x) = \dfrac{1}{\sqrt{2\gamma \pi}x}\exp\left\{ -\dfrac{(\log x-k)}{2\gamma}\right\},
\] 
with $\gamma = \gamma(\z) = \sigma^2/\mu(\z)$ and $k = k(\z) = \log x_L(\z) - \gamma(\z)$. This distribution is characterised by slim tails with exponential decay. On the contrary, von Bertalanffy-type growths, obtained from \eqref{eq:Phi} with $-1\le\delta(\z)<0$, are associated to Amoroso-type distributions
\[
f^\infty(\z,x) = \dfrac{|\delta|}{\Gamma(k/|\delta|)}\dfrac{\theta^{k}}{x^{k+1}}\exp\left\{ -\left(\dfrac{\theta}{x} \right)^{|\delta|}\right\},\qquad k(\z) = \dfrac{1}{\gamma|\delta|}+1,\quad \theta(\z) = x_L(\z)\left( \dfrac{1}{\gamma\delta^2} \right)^{1/|\delta|}, 
\]
where again $\gamma = \gamma(\z) = \sigma^2/\mu(\z)$. It is important to remark that the emerging equilibrium distribution in the case $\delta<0$ exhibits fat tails with polynomial decay. From a phenomenological point of view this is a substantial difference, since fat-tailed distributions are associated to a higher probability that the tumour is large. We point the interested reader to \cite{PTZ} for more details. 

\section{Observable effect of therapeutic protocols}
\label{sect:ther}

{In the following, we interface the natural growth mechanisms under clinical uncertainties with a superimposed therapeutical protocol that seeks to steer tumours' size towards a prescribed target. Hence, at each transitions, the tumours' size is influences by two competing dynamics, the first characterized by the uncertain growth, and the second by therapeutical protocols. }
{In details,} to determine measurable effects of therapies on growth dynamics, we include a deterministic external action as an instantaneous correction of the microscopic interaction. In details, we distinguish two types of volume updates acting on the tumour growth:
\begin{itemize}
\item[$i)$] the first is based on the transition law discussed in \eqref{eq:trans}
\item[$ii)$] the second is the therapy that acts in reducing the volume of the tumour 
\be\label{eq:ther}
x^{\prime\prime} = x + \epsilon S(x)u, 
\ee
where $u\in\mathcal{U}$, where $\mathcal U$ is the set of admissible controls such that $x^{\prime\prime}\ge 0$ and $u$ is a control defined such that
\be\label{eq:control}
u = \textrm{arg}\min_{u\in \mathcal U} J(x^{\prime\prime},u),
\ee
subject to the constraint \eqref{eq:ther}. We consider also a cost function of the form  
\be\label{eq:J_general}
J(x^{\prime\prime},u)= (x^{\prime\prime}-x_d)^2 + \nu |u|^p, 
\ee
with $\nu >0$ a penalisation coefficient and $x_d>0$ the desired tumours' size reachable with the implemented therapeutical protocol. The function $S(\cdot)$ acts selectively with respect to the tumour size. 
\end{itemize}

{In the introduced framework, we highlight that the control obtained from \eqref{eq:control} subject to \eqref{eq:ther} is independent on $\z$. Furthermore,} it is worth to remark that the typical choices for the cost function \eqref{eq:J_general} are obtained for  $p = 1,2$. More general convex functions may be considered leading often to problems that are not analytically treatable. Furthermore, in the following we will concentrate on three possible selective functions $S(x) = 1,\sqrt{x}$. 

The kinetic equation expressing the control strategy  defined in \eqref{eq:trans} and in \eqref{eq:ther} is as a sum of transition operators
\be\label{eq:kinetic_ther_strong}
\partial_t f(\z,x,t) = Q_G(f)(\z,x,t) + Q_C(\z,x,t), 
\ee
where $Q_G(\cdot)$ has been defined in \eqref{eq:QG} and the influence of therapeutical protocols on the dynamics is expressed by the new operator $Q_C(\cdot)$ whose strong formulation is given by 
\be\label{eq:QC}
Q_C(f)(\z,x,t) =  \int_{\mathbb R^+} \dfrac{1}{{}^{\prime\prime}\mathcal J} f(\z,{}^{\prime\prime}x,t)dx - f(\z,x,t), 
\ee
{with ${}^{\prime\prime}\mathcal J$ the absolute value of the Jacobian of the transformation from ${}^{\prime\prime}x$ to $x$.  }
Under suitable hypotheses is possible to obtain explicit formulation of the operator $Q_C(\cdot)$ by solving the control problem \eqref{eq:control} in feedback form at the cellular level. As before equation \eqref{eq:kinetic_ther_strong} can be fruitfully rewritten in weak form 

\be\label{eq:kinetic_ther_weak}
\dfrac{d}{dt} \int_{\mathbb R_+} \varphi(x) f(\z,x,t)dx =  \int_{\mathbb R_+} \left\langle\varphi(x^\prime)-\varphi(x) \right\rangle f(\z,x,t)dx  + \int_{\mathbb R_+} (\varphi(x^{\prime\prime})-\varphi(x))f(\z,x,t)dx.
\ee
The evolution of macroscopic quantities in the constrained setting is determined by suitable choices of the test function $\varphi$. In the following we will consider two main cases based on the minimisation of the cost \eqref{eq:J_general} with $p = 1,2$.

\subsection{The case $p=2$}
Let us consider $p=2$ in the cost function \eqref{eq:J_general}. 
The minimisation of \eqref{eq:control} can be classically done by resorting to a Lagrangian multiplier approach. We recall for related approaches the works \cite{Albi0,Albi1}. We consider the Lagrangian
\be \label{eq:lagrangian}
\mathcal{L} (u,x^{\prime\prime})=J(x^{\prime\prime},u)+\alpha [ x^{\prime\prime} - x -\epsilon S(x)u  ],
\ee
where $\alpha\in\mathbb R$ is the Lagrange multiplier associated to the constraint \eqref{eq:ther}. Hence, the optimality conditions are the following
\[
\begin{cases}\vspace{0.25cm}
\dfrac{\partial}{ \partial u} \mathcal{L}(x^{\prime\prime},u)=2\nu u - \alpha  \epsilon  S(x) = 0 \\
\dfrac{\partial}{ \partial x^{\prime\prime}}  \mathcal{L}(x^{\prime\prime},u)=  2(x^{\prime\prime}-x_d)+ \alpha = 0,  \\
\end{cases}
\]
whence we find the optimal value
\be
\label{eq:optadd}
u^* = - S(x) \dfrac{\epsilon}{\epsilon^2S^2(x)+\nu}(x-x_d).
\ee
Therefore,  plugging the optimal control \eqref{eq:optadd} defined at the cellular level into \eqref{eq:control}, we obtain the controlled transition
\[
x^{\prime\prime} = x - \dfrac{\epsilon^2 S^2(x)}{\epsilon^2 S^2(x)+\nu}(x-x_d).
\]
In this way we can study the evolution of the kinetic distribution function solution of \eqref{eq:kinetic_ther_strong}-\eqref{eq:QC} through standard methods of kinetic theory. In details, we will study the evolution of observable quantities in presence of uncertain quantities. The interplay of the introduced control with epistemic uncertainties is of paramount importance to define robust protocols.

\subsubsection{Main properties}
We define the first order moment $m(\z,t)$ and the second order moment $E(\z,t)$, or energy, respectively as 
\[
m(\z,t) = \int_{\mathbb R_+} x f(\z,x,t)dx
\]
\[
E(\z,t)=\int_{\mathbb R_+} x^2 f(\z,x,t)dx,
\]
whose evolutions are obtained by considering $\varphi(x)=x,x^2$ in \eqref{eq:kinetic_ther_weak}. 

A convenient insight on the evolution of the first order moment $m(\z,t)$ can be obtained by scaling $\nu = \epsilon \kappa$, $\kappa>0$. Under the introduced hypotheses we get
\[
\dfrac{d}{dt}m(\z,t) = \dfrac{1}{\epsilon} \left \langle \int_{\mathbb R_+} (\Phi^\epsilon_\delta(x/x_L,\z)x + x\eta_\epsilon )f(\z,x,t)dx \right\rangle - \int_{\mathbb R_+}\dfrac{ S^2(x)}{\epsilon S^2(x) + \kappa}(x-x_d)f(\z,x,t)dx. 
\]
Therefore, in the time-scale $\tau = \epsilon t$, by indicating $m(\z,\tau) = m(\z,t/\epsilon)$, we get in the limit $\epsilon \rightarrow 0^+$
\[
\dfrac{d}{d\tau}m(\z,\tau)= \int_{\mathbb R_+ }\Phi_\delta(x/x_L,\z)x f(\z,x,\tau)dx- \int_{\mathbb R_+}\dfrac{ S^2(x)}{ \kappa}(x-x_d)f(\z,x,\tau)dx. 
\]
Arguing as before for the energy $E(\z,t)$ in the case of zero diffusion, i.e., with $\eta_\epsilon \equiv 0$ in \eqref{eq:trans}, we obtain 
\[
\dfrac{d}{d\tau}E(\z,\tau)= \int_{\mathbb R_+ }\Phi_\delta(x/x_L,\z)x^2 f(\z,x,\tau)dx- \int_{\mathbb R_+}\dfrac{ S^2(x)}{ \kappa}x(x-x_d)f(\z,x,\tau)dx. 
\]

Assuming $f(\z,x,t) \in L^1(\mathbb R_+)$ it is possible to show that the model \eqref{eq:kinetic_ther_weak} has an unique equilibrium distribution $f^\infty(\z,x)$, we point the interested reader to \cite{PT} (Proposition 2.1). Hence, under the introduced regularity assumption, we can obtain some information on the large time behaviour of the first and second order moment, corresponding to the quantities $m^\infty(\z)$ and $E^\infty(\z)$. In the following we discuss the effect of the introduced control by considering different selective functions:

\begin{itemize}
\item[$a)$] if $S(x) = 1$ the asymptotic mean is solution of the following identity
\[
\int_{\mathbb R_+} \Phi_\delta(x/x_L,\z)x f^\infty(\z,t)dx = \dfrac{1}{\kappa}(m^\infty(\z)-x_d).
\]
Note that since $\Phi_\delta$ is bounded for all $\z \in \mathbb R^d$ by the following uncertain quantities
\[
-\dfrac{\mu}{1+\lambda}\le \Phi_\delta\le \dfrac{\mu}{1-\lambda}
\]
we have
\be\label{eq:modulo}
\left|\int_{\mathbb R_+} \Phi_\delta(x/x_L,\z)x f^\infty(\z,t)dx \right| \le \int_{\mathbb R_+} \left| \Phi_\delta(x/x_L,\z) \right| x f^\infty(\z,t)dx \le \dfrac{\mu}{1-\lambda}m^\infty(\z).
\ee
Hence, the following inequality holds
\[
\dfrac{1}{\kappa} |m^\infty(\z) - x_d| \le \dfrac{\mu}{1-\lambda}m^\infty(\z),
\]
whose solution is such that
\[
\dfrac{1-\lambda}{1-\lambda+\kappa\mu} x_d \le m^\infty(\z) \le \dfrac{1-\lambda}{1-\lambda-\kappa\mu}x_d
\]
provided $\kappa < \min_{\z \in \mathbb R^d} \frac{1-\lambda}{\mu}$. We easily observe that in the limit $\kappa\rightarrow 0^+$ corresponding to vanishing penalisation of the control the large time mean size is such that $m^\infty(\z)\rightarrow x_d$. In other words, we have
\[
-\dfrac{\kappa\mu}{1-\lambda+\kappa\mu}x_d \le m^\infty(\z) - x_d \le \dfrac{\kappa\mu}{1-\lambda-\kappa\mu}x_d,
\]
and
\begin{equation}\label{eq:est}
|m^\infty (\z)- x_d| \le  \dfrac{\kappa\mu}{1-\lambda-\kappa\mu}x_d.
\end{equation}

Let us assume that $\sigma^2 = 0$. Then, the second order {moment} is such that 
\[
\int_{\mathbb R_+} \Phi_\delta(x/x_L,\z)x^2 f^\infty(\z,t)dx = \dfrac{1}{\kappa}(E^\infty(\z)-m^\infty(\z)x_d).
\]
We note that
\[
\left|\int_{\mathbb R_+} \Phi_\delta(x/x_L,\z)x^2 f^\infty(\z,t)dx \right| \le \int_{\mathbb R_+} \left| \Phi_\delta(x/x_L,\z) \right| x^2 f^\infty(\z,t)dx \le \dfrac{\mu}{1-\lambda}E^\infty(\z),
\]
since $\Phi_\delta$ is bounded for all $\z \in \mathbb R^d$, as we observed before. Consequently, we have 
\[
E^\infty(\z)-m^\infty(\z)x_d \leq \dfrac{\kappa\mu}{1-\lambda}E^\infty(\z),
\]
that is
\[
E^\infty(\z)\left( 1 - \dfrac{\kappa\mu}{1-\lambda} \right) \leq m^\infty(\z)x_d.
\]
Since in the limit $\kappa\rightarrow0^+$ we have observed that $m^\infty(\z)\rightarrow x_d$ we can write
\[
E^\infty(\z) - (m^\infty(\z) )^2 \leq 0.
\]
We observe also that $E^\infty(\z) - (m^\infty(\z) )^2\geq0$ by definition, since it is the variance of the random variable $X\sim f^\infty(\z,x)$. In other words, in the limit $\kappa\rightarrow0^+$ we have $\textrm{Var}_{f^\infty}[X] \rightarrow 0$, that is, the equilibrium distribution tends to a Dirac delta centred in $x=x_d$.

\item[$b)$] We consider now the case $S(x) = \sqrt{x}$ corresponding to a heavier control on large sized tumours. We can observe that in this case the asymptotic first order {moment} solves 
\[
\int_{\mathbb R_+} \Phi_\delta(x/x_L,\z) x f^\infty(\z,x)dx = \dfrac{1}{\kappa} \int_{\mathbb R_+}x(x-x_d)f^\infty(\z,x)dx. 
\]
In details, since from the Jensen's inequality we have
\[
\int_{\mathbb R_+}(x-x_d)^2 f^\infty(\z,x) dx \ge \left(\int_{\mathbb R_+}(x-x_d)f^\infty(\z,x)dx\right)^2. 
\]
we get
\[
\dfrac{1}{\kappa}m^\infty(\z)(m^\infty(\z)-x_d) \le \int_{\mathbb R_+}\Phi_\delta(x/x_L,\z)xf^\infty(\z,x)dx. 
\]
Therefore, thanks to \eqref{eq:modulo} we obtain 
\begin{equation}\label{eq:est_b}
|m^\infty(\z) - x_d| \le \dfrac{\mu\kappa}{1-\lambda}. 
\end{equation}
As obtained in point $(a)$ we obtain that for vanishing penalisation $\kappa\rightarrow 0^+$ the asymptotic first order {moment} is such that $m^\infty(\z) \rightarrow x_d$.

Assuming now that $\sigma^2 = 0$ the asymptotic energy solves  
\[
\int_{\mathbb R_+} \Phi_\delta(x/x_L,\z) x^2 f^\infty(\z,x)dx = \dfrac{1}{\kappa} \int_{\mathbb R_+}x^2(x-x_d)f^\infty(\z,x)dx,
\]
from which we get
\[
\left| (m^\infty(\z))^3 - x_d E^\infty(\z) \right| \le \dfrac{\mu\kappa}{1-\lambda}E^\infty(\z),
\]
and in the limit $\kappa\rightarrow 0^+$ we obtain that the large time distribution tends to a Dirac delta centred in $x=x_d$.
\end{itemize}

In both the discussed cases and in particular from \eqref{eq:est} and \eqref{eq:est_b}, we can observe that the introduced protocols induce the mean tumours' sizes to stick the deterministic target size $x_d$. These results have an important consequence on the uncertainties of the system. In particular, looking at the variance with respect to $\z \in \mathbb R^d$ we have
\[
\textrm{Var}_\z(m^\infty(\z)) = \textrm{Var}_\z(m^\infty(\z)-x_d) = \mathbb E_\z[(m^\infty(\z)-x_d)^2] -\mathbb E_\z[(m^\infty(\z)-x_d)]^2,
\]
from which we get
\begin{equation}
\textrm{Var}_\z(m^\infty(\z)) \le \mathbb E_\z[(m^\infty(\z)-x_d)^2] \le \textrm{max}\left\{\dfrac{\kappa\mu}{1-\lambda-\kappa\mu}x_d,\dfrac{\mu\kappa}{1-\lambda} \right\}. 
\end{equation}
Hence, since $\textrm{max}\left\{\dfrac{\kappa\mu}{1-\lambda-\kappa\mu}x_d,\dfrac{\mu\kappa}{1-\lambda} \right\} \rightarrow 0$ for $\kappa\rightarrow 0^+$, we argue that the introduced controls are capable to dampen invariably the variability due to the presence of clinical uncertainties $\z \in \mathbb R^d$.

\subsubsection{Large time behaviour of the controlled model} \label{subsec:FPA}

At this point, proceeding as in Section \ref{subsec:equilibria} for the new kinetic model \eqref{eq:kinetic_ther_strong} we can assess the effects of the control therapies on the emerging kinetic distribution. In the limit $\epsilon \rightarrow 0^+$ and scaling $\nu = \epsilon \kappa$, where $\kappa>0$ is the scaled penalisation, the kinetic equation converges to a Fokker-Planck equation with modified drift term that takes into account the presence of the control. The resulting Fokker-Planck-type equation reads
\be \label{eq:FPcontr}
\partial_t f(\z,x,t) = \partial_x \left[ -\Phi_\delta(x/x_L,\z)xf(\z,x,t) + \dfrac{\sigma^2}{2}\partial_x(x^2 f(\z,x,t)) \right] + \dfrac{1}{\kappa} \partial_x \left[S^2(x)(x-x_d) f(\z,x,t) \right]. 
\ee

Since we have obtained in Section \ref{subsec:equilibria} that if $\delta<0$ the introduced model {leads} to equilibrium distributions with polynomial tails, linked to a high probability that the tumours' sizes are large, we concentrate on this case. Under this assumption, the asymptotic large time distribution of the controlled model is given by
\[
f^\infty(\z,x) = C_{\mu,\sigma^2,x_L}(\z) \left(\dfrac{1}{x}\right)^{\frac{1}{\gamma|\delta|}+2} \textrm{exp}\left\{ - \dfrac{2}{\sigma^2\delta^2} \left(\dfrac{x}{x_L} \right)^\delta \right\} \textrm{exp}\left\{-\dfrac{2}{\sigma^2\kappa}\int \dfrac{S^2(x)(x-x_d)}{x^2}dx \right\}, 
\]
with $C_{\mu,\sigma^2,x_L} >0$ a normalisation constant. 
Hence, if $S(x) = 1$ a direct computation gives
\[
f^\infty(\z,x) = C_{\mu,\sigma^2,x_L}(\z) \left( \dfrac{1}{x}\right)^{\frac{1}{\gamma|\delta|}+\frac{2}{\sigma^2\kappa}+2}\textrm{exp}\left\{ - \dfrac{2}{\sigma^2\delta^2} \left(\dfrac{x}{x_L} \right)^\delta \right\} \textrm{exp}\left\{ -\dfrac{2x_d}{\sigma^2\kappa}\dfrac{1}{x} \right\}, 
\]
and the emerging equilibrium of the controlled exhibits again power law tails for large $x$'s. Anyway, it is worth to observe that the exponent increases due to the presence of the introduced control.  On the other hand, for selective  controls with $S(x) = \sqrt{x}$ we get 
\[
f^\infty(\z,x) = C_{\mu,\sigma^2,x_L}(\z) \left(\dfrac{1}{x}\right)^{\frac{1}{\gamma|\delta|}-\frac{2x_d}{\sigma^2\kappa}+2}\textrm{exp}\left\{ - \dfrac{2}{\sigma^2\delta^2} \left(\dfrac{x}{x_L} \right)^\delta \right\} \textrm{exp}\left\{ -\dfrac{2}{\sigma^2\kappa}x \right\},
\]
provided $\kappa>2x_d\gamma|\delta|/\sigma^2$, corresponding to a distribution with exponential decay of the tail.  In other words, even if the introduced therapies are capable to reduce in any case the influence of clinical uncertainties, selective-type controls, whose action is heavier on large tumours, are necessary to modify the nature of the emerging distribution of tumours' sizes.  

\subsection{The case $p=1$}
Let us consider $p=1$ in \eqref{eq:J_general}. Proceeding as before we consider the Lagrangian \eqref{eq:lagrangian} with the cost function $J(x^{\prime\prime},u)= (x^{\prime\prime}-x_d)^2 + \nu |u|$. The optimality conditions read now
\[
\begin{cases}\vspace{0.25cm}
\dfrac{\partial}{ \partial u} \mathcal{L}(x^{\prime\prime},u)=\nu \sign(u) - \alpha  \epsilon  S(x) = 0 \\
\dfrac{\partial}{ \partial x^{\prime\prime}}  \mathcal{L}(x^{\prime\prime},u)=  2(x^{\prime\prime}-x_d)+ \alpha = 0.  \\
\end{cases}
\]
A direct solution of the previous system leads, as before, to a feedback formulation of the optimal control that can be written as follows
\[
u^*=\Pi_{U}\left(\mathbb{S}^1(x-x_d)\right),
\]
where the operator $\mathbb{S}^1(x-x_d)$ is defined as
\[
\mathbb{S}^1(x-x_d) :=
\begin{cases}\vspace{0.25cm}
\dfrac{x-x_d}{|x-x_d|} \dfrac{\nu}{2\epsilon^2 S^2(x)} - \dfrac{x-x_d}{\epsilon S(x)}, \quad\qquad & |x-x_d| > \dfrac{\nu}{2\epsilon S(x)}, \\
0 & \textrm{otherwise} \\
\end{cases}
\]
and $\Pi_{U}$ is the projection onto a compact subset $U\subset\mathbb{R}$. It should be noted that for any value of $\epsilon>0$ the applied control is active only on a portion of tumours. This result is coherent with analogous works in related fields, see \cite{AFK,DPT}. 

In this scenario, it is not possible to apply the same arguments of Section \ref{subsec:FPA} to get analytical results on the evolution of observable quantities. Furthermore, the derivation of surrogate  Fokker-Planck-type models does not help to obtain insights on the large time behaviour of the system. As a consequence, in the following we will focus on the consistent numerical approaches to have a qualitative indication of the emerging phenomena.

\section{Quantities of interest and data} \label{sec:QoI}

In this section, we face the calibration of the kinetic model \eqref{eq:kin_strong} defined in Section \ref{subsec:equilibria} in presence of uncertain quantities by means of experimental data. In particular, to obtain some evidence on the distribution of uncertain quantities, we focus on the microscopic laws defined in \eqref{eq:trans} to get a patient-based estimation of all the relevant parameters characterising the dynamics. Thus, to deal with the uncertainties brought by the parameter $\z$ and affecting the evolution of the distribution $f(\z,x,t)$, we consider as a quantity of interest (QoI) the expected evolution of the first order {moment} $\mathbb E_\z[m(\z,t)]$. In this way, we are able to compare the theoretical and numerical results with the measures of our dataset relative to the cohort of subjects affected by glioblastoma. Empirical measurements of a subject's tumour sizes correspond to a specific realisation of a particular value of the random variable $\z$. Therefore, the average behaviour of a glioblastoma is the result of the superposition of different dynamics, produced by different values of $\z$, incorporating the subjects' variability, that are then weighted by the associated probability measure $\rho(\z)$.

In particular, we are interested in the analytical and numerical solutions obtained for $\delta  \rightarrow 0$ and for $\delta <0$, reproducing Gompertz and von Bertalanffy growth models respectively. Parameter estimation in tumour growth dynamics is a classical problem and we mention \cite{Laird, Laird1} for an introduction on the topic. More recently a similar problem has been considered for glioblastoma in \cite{Za}. 

\subsection{Dataset construction and Segmentation}

In this work, we consider clinical data for tumour growth relative to a cohort of patients referred to IRCCS Mondino, collected from 2011 to 2021. Among 263 subjects suffering from brain tumour, we select those affected by primary glioblastoma. In all these cases MRI (Magnetic Resonance Imaging) scans were available after each visit. Combining the neuro-radiological and the clinical information, we choose among the selected patients the ones who exhibit an initial tumour free-growth and that have at least two MRI scans at different times. In this way, we are capable to estimate the patient-based growth dynamics. Anyway, only very few observations can be obtained of these characteristics since the great majority of patients are enrolled for follow-up at Mondino after initial treatments. For these reasons, we include subjects with treatments' interruption. At the end of this preliminary analysis, we considered the evolution of the tumours' size of 13 patients.

Among the patients' MRI sequences, typical of the MRI brain tumour acquisition protocol,  we are interested in the T1 weighted 3D MRI scans with contrast agent to estimate the subjects' tumour volume $x$ in $mm^3$ at a given time. The T1-weighted MRI images rely upon longitudinal relaxation of the tissue's magnetisation vector due to the protons spin-lattice interactions. Different tissue types are characterised by different T1 relaxation times, therefore it is possible to differentiate anatomical structures. An injection of a contrast agent, such as gadolinium, during the T1- weighted image acquisition, supplies information about current disease activity. In fact, passing through the blood brain barrier, the contrast agent reveals inflammation areas that appear brighter, helping in identifying the tumours' contours.

For each subject, the glioblastoma volume segmentation is performed using the software 3D slicer \cite{3Dslicer}. We combine a data clustering algorithm and manual segmentation corrections. In particular, we apply the region growing algorithm based on the examination of neighbouring pixels of the initial seeds, a set of selected points in the region of interest, determining whether a neighbour pixel should be added to that region or not. After that, a manual correction of contours is performed. The procedure is iterated in the axial, coronal and sagittal image projection in order to obtain more precise results. To determine the tumours' volume, the number of voxels contained in the segmentation and the MRI metadata information have been considered.

\subsection{Growth curves and growth model parameters' distributions} 
\label{subsec:growth}
To determine the empirical distribution of the parameters characterising the tumours' dynamics we adopted a two-level approach. {In particular, in the free-growth scenario, we estimated the parameters of \eqref{eq:trans}.} This estimation will be then kept to evaluate the observable effects of the treatment. In more details, for a cohort of $N$ patients we define $\{\hat x_i(t^n)\}_{i=1}^N$ the observed volume size at time $t^n$. 

Assuming Gompertz-type growths we need to estimate a 2D vector for each patient, i.e. the tumour growth rate $\alpha>0$ and the carrying capacity $x_L>0$. We indicate with $\Theta = (\alpha,x_L)$ the 2D vectors of parameters. Hence, in the time interval $[0,T]$ we solve a least square problem based on the minimisation of a suitable norm of theoretical and empirical tumour's sizes measured at the available times $t^0,\dots,t^n\le T$. More precisely, we considered a minimisation problem based on the following norm
\begin{equation}
\label{eq:norm_data}
\min_{\Theta} \left[{ \sum_{h \in H_i}| x_i(t^h)-\hat x_i(t^h)|} + \beta \| \Theta\|_{L^1} \right],
\end{equation}
{where $H_i$ collects all the observations of the tumour's volume of the $i$th subject}. Furthermore, we introduced the regularisation parameter $\beta>0$. In the case $\delta\rightarrow 0^+$ we considered the theoretical evolution for $x$ given by \eqref{eq:gomp}.

For von Bertalanffy-type dynamics we have to estimate a 3D vector for each patient $\Theta = (a,p,q)$, with $a = \delta+1$, as observed in Section \ref{sect:transition}. Furthermore, information on the carrying capacity $x_L$ has been considered compatible with the Gompertz case. In the time interval $[0,T]$ we solved a least square problem \eqref{eq:norm_data} where the theoretical evolutions of the tumours' volumes are given by \eqref{eq:vb}.

\begin{figure}
	\centering
	\includegraphics[scale = 0.50]{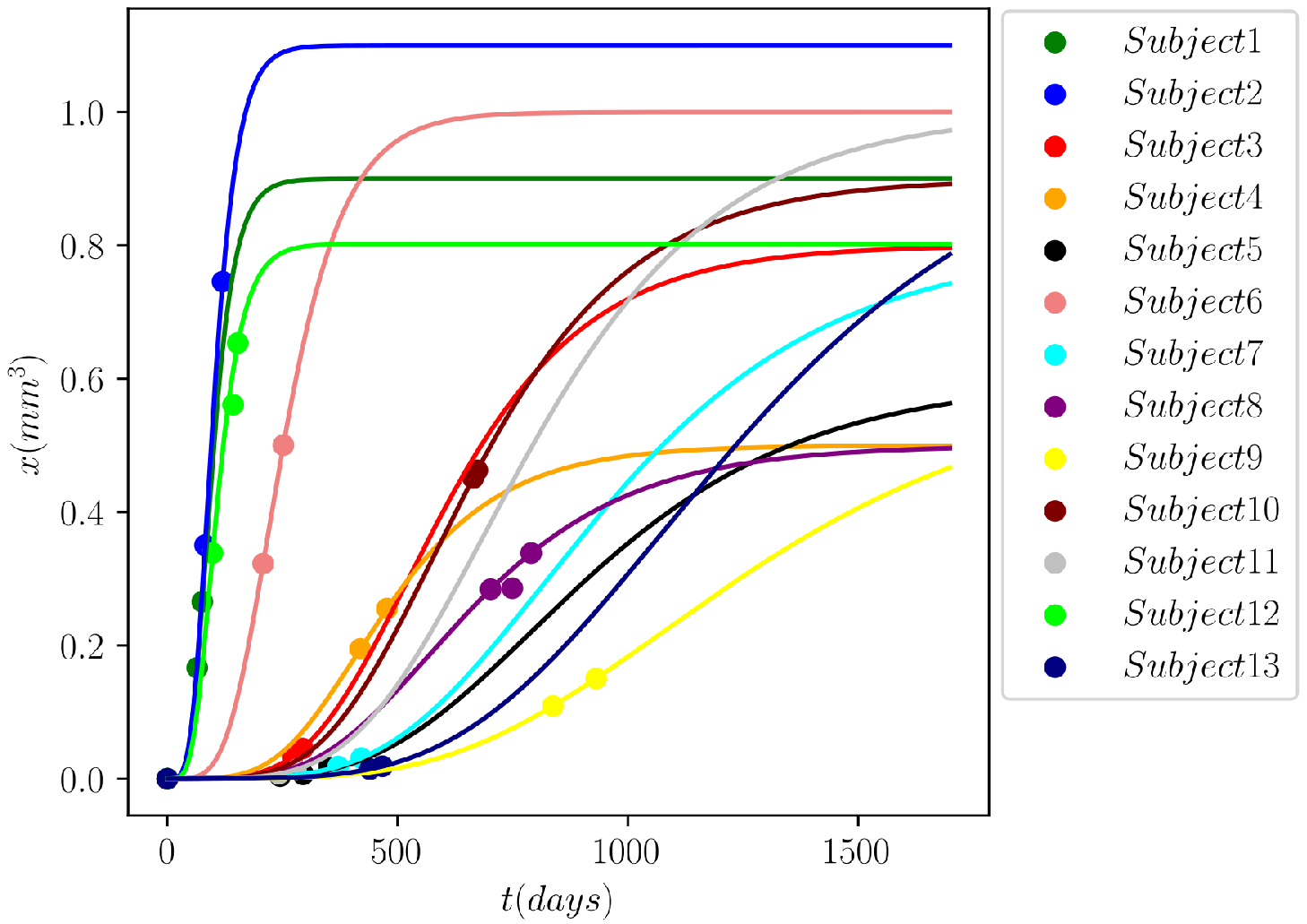}
	\includegraphics[scale = 0.50]{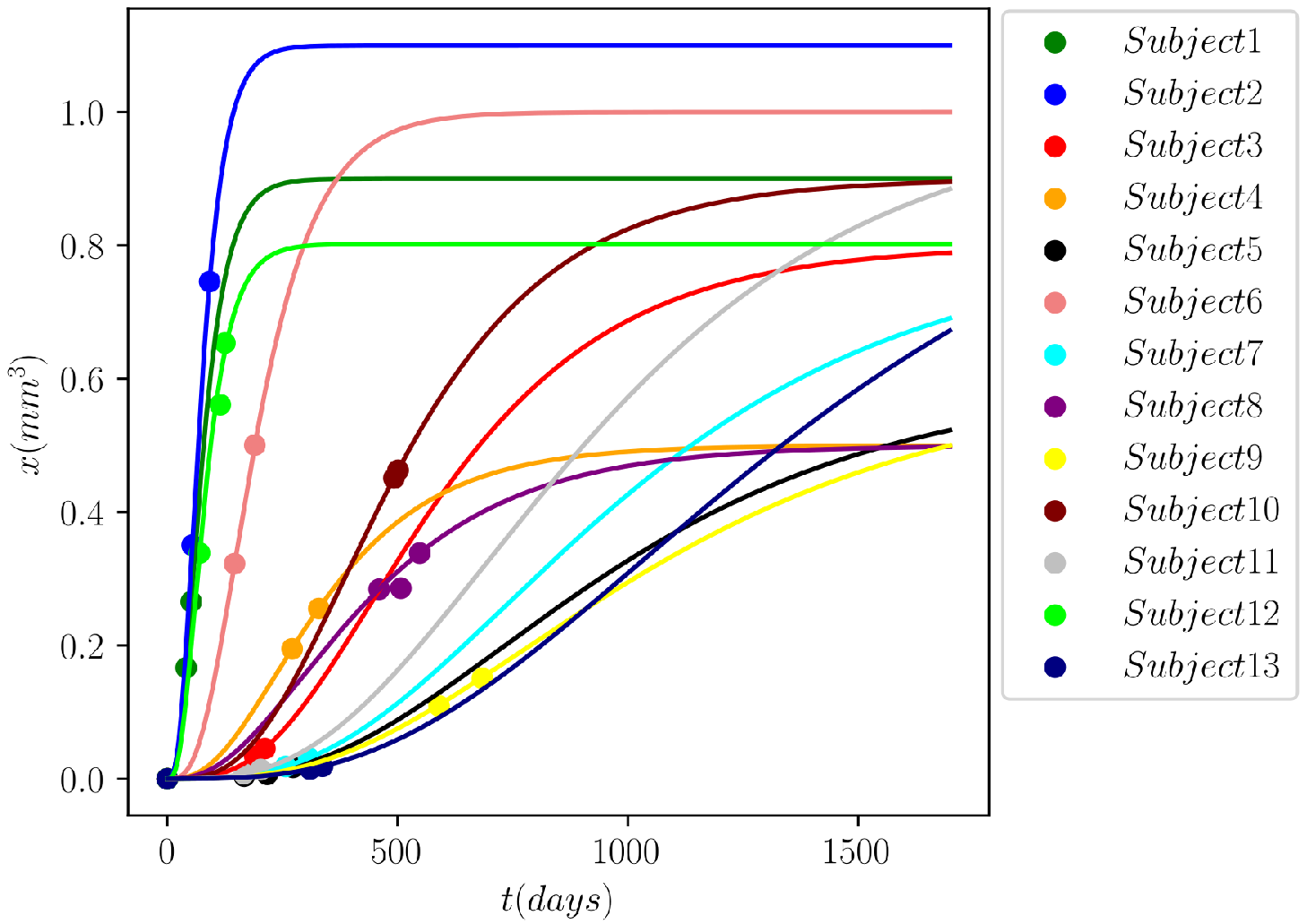}
		\caption{Growth curves and empirical data relative to untreated glioblastoma. The trajectory of each curve (solid line) and the empirical volume size (circle marker) data are shown for each patient and for different growth laws: Gompertz case (left plot), von Bertalanffy case (right plot).
		Values reported on y-axes are scaled by a quantity of $10^5$.}
	\label{fig:curve_sperimentali}
	
\end{figure}

Since the first MRI time point $t_0$ and tumour size $x_0$ are different for each subject, we need to find a common point with the aim of comparing the patient specific growth curves for both models. As initial volume, we take the tumour size $1mm^3$ as the mentioned point. This choice is justified by the fact that the smallest appreciable MRI voxel dimension is $1mm^3$. Hence, we solve through standard numerical methods \ref{eq:micro} to obtain, for each subject, the specific time corresponding to $1mm^3$. Subsequently, we translate for each subject the initial time of the estimated time. The obtained growth curves and empirical volume size data are shown for each subject and for Gompertz and von Bertalanffy models in Figure \ref{fig:curve_sperimentali}.

\begin{figure}
	\centering
	\includegraphics[scale = 0.3250]{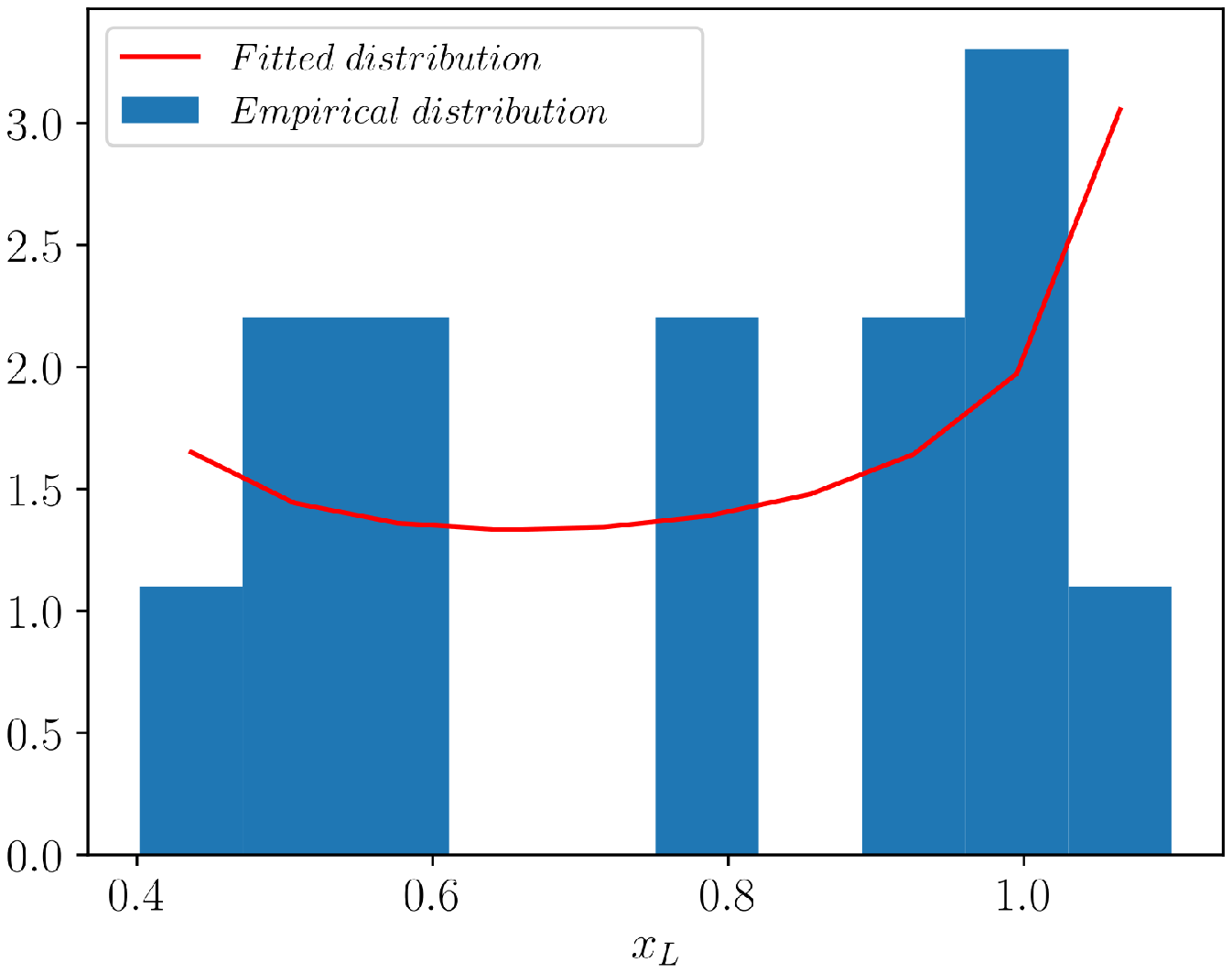}
	\includegraphics[scale = 0.3250]{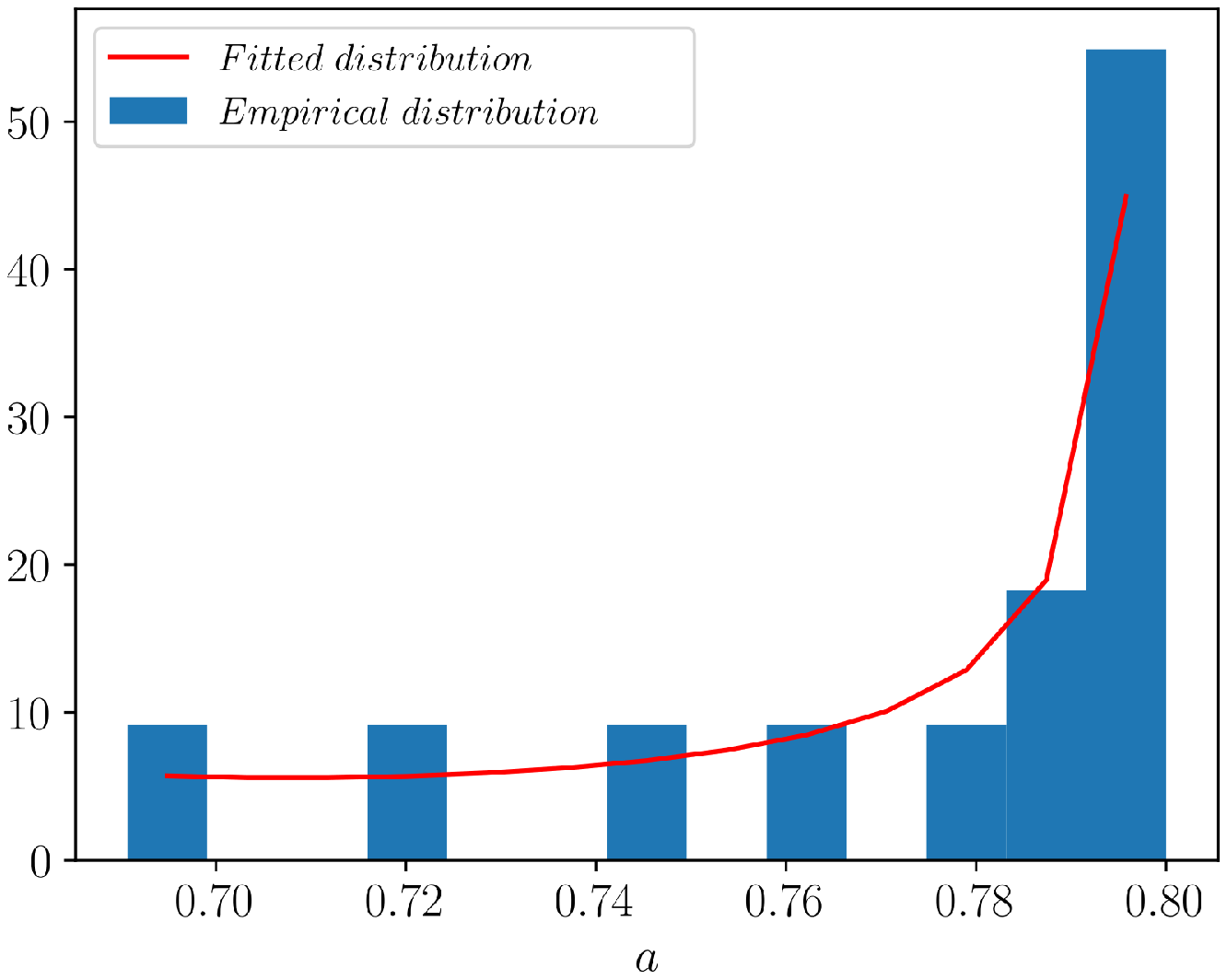}
	\includegraphics[scale = 0.3250]{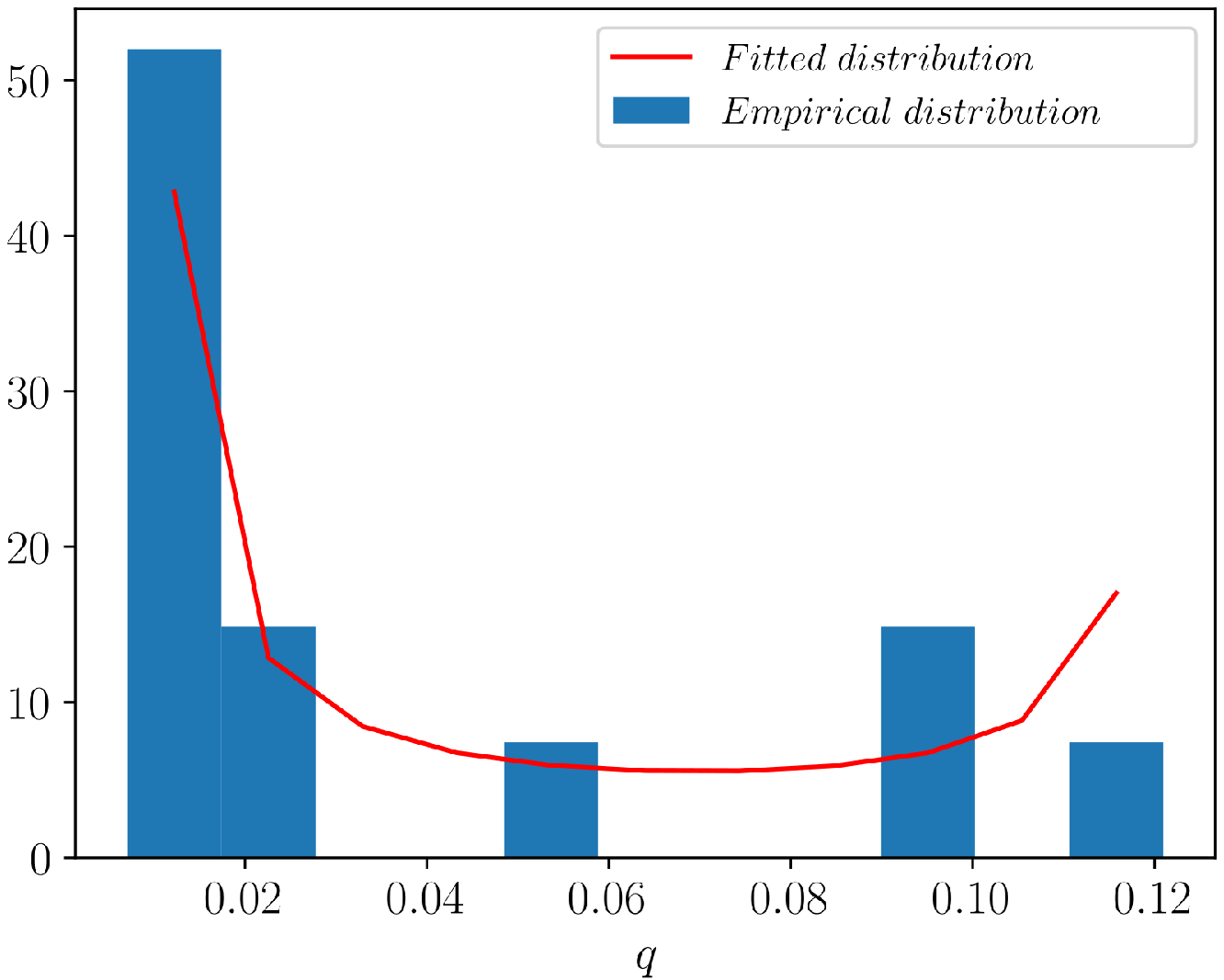}
	\caption{Empirical distributions of the obtained parameters and fitted Beta distributions with parameters given in Table  \ref{table:tabella}.}
	\label{fig:distribuzioni}
\end{figure}

To understand the trends of the aforementioned model parameters, incorporated in the random variable $\z$, we construct the associated histograms and we determine the theoretical distributions that better reproduce each of them in the associated range of variability. The results are shown in the Figures \ref{fig:distribuzioni}. We obtained a poor fit of the parameter $\alpha$ characterising Gompertz-type growths and we decided to consider an uninformative uniform distribution over the observed interval of variability $[0.001,0.03]$. Anyway, we observe that the range of $\alpha$ is consistent with values reported in \cite{Za} and obtained from a global fit on a larger subjects data cohort.

The parameters of the theoretical distributions are obtained by maximising the proper likelihood function. To verify the goodness of the theoretical representations, we quantify the distance between each distribution function of the empirical sample and the cumulative distribution function of the selected theoretical one by performing the Kolmogorov-Smirnov test. The results are summarised in Table \ref{table:tabella}.

\begin{table}
	\centering
\begin{tabular}{c|c|c|c|c}
	Parameter & Range & Distribution & Constants $(c_1,c_2)$ & KS-pvalue\\
	\hline
	$x_L$ & $[0.4,  \,1.1]$  & Beta &  $(0.705,\,0.574)$ & $0.823$\\
	$a$   & $[0.69, \,0.8]$  & Beta &  $(0.656,\,0.193)$ & $0.902$ \\
	$q$   & $[0.007,\,0.12]$ & Beta &  $(0.112,\,0.238)$ & $0.314$
\end{tabular}
\caption{We report for each parameter the best fitted Beta distribution characterised by the constants $(c_1,c_2)$ in the third column and the range of definition in the second column (the $x_L$ range is scaled by a quantity of $10^5$). The quantification of the goodness of the theoretical representations is given by KS-pvalue in the last column.}
\label{table:tabella}
\end{table}

All other parameters, as can be seen from Table \ref{table:tabella}, are instead well described by a Beta distribution defined by 
\[
f(x,c_1,c_2)= \frac{\Gamma(c_1+c_2)x^{c_1-1} (1-x)^{c_2-1}}{\Gamma(c_1)\Gamma(c_2)}
\]
with $c_1$ and $c_2$ the shape parameters that have been reported in the third column of Table \ref{table:tabella}.

\section{Numerical tests}\label{sect:numerics}
In this section we introduce accurate numerical strategies for kinetic equations based on a stochastic Galerkin formulation of the derived equations, see \cite{CZ,DPZ,Par,ZJ} and the references therein. In particular, we present several numerical tests highlighting the obtained theoretical results focusing first on the untreated tumour growth case providing results on spectral convergence of the adopted methods.  Furthermore, we compare the evolution of the QoI with the experimental data. Next, we investigate the case including therapies through the considered control protocols testing its effectiveness in damping the introduced uncertainties at the level of observable quantities. In the following, we will consider all the tumours' volumes scaled by a factor of $10^5$. 

\subsection{Stochastic Galerkin methods}

The stochastic Galerkin (sG) method is based on the construction of a set of hierarchical orthogonal polynomials $\{\Psi_k(\z)\}^{M}_{k=0}$ of degree less or equal to $M\in \mathbb N$, orthonormal with respect to the PDF of the random parameters $ \rho(\z)$, that is
\[
\mathbb E_{\z}[\Psi_k(\z)\Psi_h(\z)]=\int_{\mathbb R^d}\Psi_k(\z)\Psi_h(\z)\rho(\z)d\z=\delta_{kh},\;\;\;\;\;k,h=0,\dots,M,
\]
where $\delta_{kh}$ is the Kronecker delta. The choice for the orthogonal polynomials obviously depends on the PDF of the parameters $\rho(\z)$ and follows the so-called Wiener-Askey scheme, see \cite{Xiu,WienAsk}.

Let $f = f(\z,x,t)$ be the solution of a Fokker-Planck equation at time $t\ge0$, provided that it is sufficiently regular, it can be approximated by $f^M$ that is defined as follows
\be
\label{eq:fapprox}
f(\z,x,t) \approx f^M(\z,x,t)=\sum_{k=0}^M\hat{f}_k(x,t)\Psi_{k}(\z),
\ee
where $\hat{f}_k(x,t)$ is the projection of the solution over the space generated by the polynomial of degree $k = 0,\dots,M$
\[
\hat{f}_k(x,t) := \mathbb E_{\z}[f(\z,x,t)\Psi_{k}(\z)] = \int_{\mathbb R^d}f(\z,x,t)\Psi_{k}(\z)\rho(\z)d\z.
\]
Hence, if we substitute the approximation \eqref{eq:fapprox} of the PDF into the Fokker-Planck equation \eqref{eq:FP}, exploiting the orthonormality of the polynomials, we find a system of $M+1$ equations for the time evolution of the projections $\hat{f}_k(x,t)$, that reads
\be
\label{eq:FPsG}
\partial_{t}\hat{f}_k(x,t) = \partial_{x}\left[ \sum_{h=0}^M\mathcal A_{kh}(x) \hat{f}_h(x,t) + \dfrac{\sigma^2}{2}\partial_{x}(x^2\hat{f}_k(x,t))\right],
\ee
where the matrix $\mathcal A_{kh}(x)$ is defined as
\[
\mathcal A_{kh}(x) = - \int_{\mathbb R^d} x \Phi_\delta(x/x_L,\z)\Psi_{k}(\z)\Psi_{h}(\z)\rho(\z)d\z.
\]
We stress the fact that the system of equations \eqref{eq:FPsG} is deterministic since it does not depend on the random parameters $\z$. The main advantage of the stochastic Galerkin approach relies on the fact that, if the solution of the PDE of interest is sufficiently regular, the approximated solution spectrally converges to the correct solution of the problem. This translates into the fact that it is sufficient to consider $M$ relatively small. 

Analogous computations can be performed in the model that includes the introduced control \eqref{eq:FPcontr} with the only difference that the drift coefficient results modified by an additional term. In particular, in the controlled case the matrix $\mathcal A_{kh}(x)$ reads
\benn
\begin{split}
\mathcal A_{kh}(x) &= - \int_{\mathbb R^d}\left( x \Phi(x/x_L,\z) - \frac{S^2(x)(x-x_d)}{\kappa}\right)\Psi_{k}(\z)\Psi_{h}(\z)\rho(\z)d\z 
\end{split}.
\eenn

In order to prove the stability result for the sG scheme, we may reformulate the Fokker-Planck equation \eqref{eq:FPsG} in a more compact form. If we define the $M+1$ vector $\hat{\f}(x,t)=(\hat{f}_0(x,t),\dots,\hat{f}_M(x,t))$, the $(M+1)\times(M+1)$ matrix $\mathbf{A}(x)=\{ \mathcal A_{kh}(x) + \sigma^2 x \mathbb I \}_{k,h=0}^M$, being $\mathbb I$ a unitary matrix,  and the diffusion coefficient $D(x)=x^2\sigma^2/2$, we have
\be \label{eq:FPvett}
\partial_t \hat{\f} (x,t) = \partial_x \left[ \mathbf{A}(x) \hat{\f}(x,t) + D(x) \partial_x  \hat{\f}(x,t)  \right].
\ee
We denote with $\| \hat{\f} \|_{L^2}$ the standard $L^2$ norm of the vector $\hat{\f}(x,t)$
\[
\| \hat{\f} \|_{L^2} := \left[ \int_{\mathbb R_+} \left( \sum_{k=0}^M \hat{f}_k^2(x,t)  \right)^2 dx \right]^{1/2},
\]
and we observe that, thanks to the orthonormality of the polynomials $\{ \Psi_k \}_{k=0}^M$ in $L^2(\Omega)$, we have 
\[
\| f^M \|_{L^2(\Omega)} = \| \hat{\f} \|_{L^2}.
\]
Now, we can show the stability result. 
\begin{theorem}
Assume that there exists two constants $C_A > 0 $ such that $\| \partial_x \mathcal A_{kh} \|_{L^\infty}\leq C_A$ for every $k,h=0,\dots,M$ and $D(x)>0 $ for every $x\in\mathbb R_+$, then 
\[
\| \hat{\f} \|_{L^2}^2 \leq e^{C_A\, t} \| \hat{\f}(0) \|_{L^2}^2. 
\]
\begin{proof}
We multiply every component of \eqref{eq:FPvett} by $\hat{f}_k$ and we integrate over $\mathbb R_+$ to get
\[
\int_{\mathbb R_+} \dfrac{1}{2} \partial_t \left( \hat{f}_k^2 \right) dx = \int_{\mathbb R_+} \hat{f}_k \partial_x \left[  \sum_{h=0}^M \mathcal A_{kh} \hat{f}_h + D(x) \partial_x \hat{f}_k \right] dx.
\]
We integrate by parts the transport term on the right-hand side of the equation to obtain
\benn
\begin{split}
   &\sum_{h=0}^M \int_{\mathbb R_+} \hat{f}_k \partial_x \left( \mathcal A_{kh} \hat{f}_h \right) dx = \sum_{h=0}^M \int_{\mathbb R_+} \left( \hat{f}_k \hat{f}_h \partial_x \mathcal A_{kh} +  \hat{f}_k \mathcal A_{kh} \partial_x \hat{f}_h  \right) dx \\
   & = - \sum_{h=0}^M \int_{\mathbb R_+} \mathcal A_{kh} \partial_x \left( \hat{f}_k \hat{f}_h \right) dx - \sum_{h=0}^M \int_{\mathbb R_+} \hat{f}_h  \partial_x \left( \mathcal A_{kh} \hat{f}_k  \right) dx. 
\end{split}
\eenn
We sum over $k=0,\dots,M$ and we exploit the symmetry of $\mathbf{A}$ to have
\benn
\begin{split}
   &2 \sum_{k,h=0}^M \int_{\mathbb R_+} \hat{f}_k \partial_x \left( \mathcal A_{kh} \hat{f}_h \right) dx =  - \sum_{k,h=0}^M \int_{\mathbb R_+} \mathcal A_{kh} \partial_x \left( \hat{f}_k \hat{f}_h \right) dx \\
   &= \sum_{k,h=0}^M \int_{\mathbb R_+} \hat{f}_k \hat{f}_h  \partial_x \mathcal A_{kh} dx.
\end{split}
\eenn
Since $\| \partial_x \mathcal A_{kh} \|_{L^\infty}\leq C_A$ and from Cauchy-Schwartz inequality, we have
\[
\sum_{k,h=0}^M \int_{\mathbb R_+} \hat{f}_k \partial_x \left( \mathcal A_{kh} \hat{f}_h \right) dx \leq \dfrac{C_A}{2} \parallel \hat{\f} \parallel_{L^2}^2.
\]
As for the diffusion term, we have 
\benn
\begin{split}
\sum_{k=0}^M\int_{\mathbb R_+} \hat{f}_k \partial_x \left( D(x) \partial_x \hat{f}_k \right) dx &= -\sum_{k=0}^M\int_{\mathbb R_+} D(x) \left(  \partial_x \hat{f}_k \right)^2 dx \le 0, 
\end{split}
\eenn 
since $D(x)\ge 0$ by assumption. 
If we sum over $k$, the left-hand side is nothing but the derivative in time of the $L^2$ norm of $\hat{\f}$
\[
\sum_{k=0}^M\int_{\mathbb R_+} \dfrac{1}{2} \partial_t \left( \hat{f}_k^2 \right) dx = \dfrac{1}{2} \partial_t \| \hat{\f} \|_{L^2}^2.
\]
Finally, we have 
\[
\dfrac{1}{2} \partial_t \| \hat{\f} \|_{L^2}^2 \leq \dfrac{C_A}{2}  \| \hat{\f} \|_{L^2}^2,
\]
and thanks to Gronwall's Lemma we conclude. 
\end{proof}
\end{theorem}

We concentrate on the case where the evolution of the Fokker-Planck equations \eqref{eq:equilibria} and \eqref{eq:FPcontr} is affected by an uncorrelated 2D random term in the Gompertz case, i.e. $\z = (z_1,z_2)$ and $\z \sim \rho(z_1,z_2) = \rho_1(z_1)\rho_2(z_2)$, or by an uncorrelated 3D random term in the von Bertalanffy case, i.e. $\z = (z_1,z_2,z_3)$ and $\z \sim \rho(z_1,z_2,z_3) = \rho_1(z_1)\rho_2(z_2)\rho_3(z_3)$. The distribution of the components of the random vectors are determined by the analysis presented in Section \ref{sec:QoI}. 

In the limit $\delta \rightarrow 0^+$, corresponding to a kinetic Gompertz model, the approximated solution is therefore given by 
\[
f(\z,x,t) = f^M(\z,x,t) \approx  \sum_{h,k=0}^M \hat f_{hk}(x,t)\Psi^1_h(z_1)\Psi^2_k(z_2)
\]
being $\{\Psi^1_h\}_{h=0}^M$ and $\{\Psi^2_k(z_2)\}_{k=0}^M$ the set of polynomials orthonormal with respect to $\rho_1(z_1)$ and $\rho_2(z_2)$ respectively. Under the introduced assumptions we obtain the following set of equations
\begin{equation}\label{eq:sG_Gompertz}
    \partial_t \hat f_{hk} = \partial_x \left[ \sum_{\ell,r = 0}^{M} \mathcal A_{hk\ell r} \hat{f}_{\ell r}(x,t) + \dfrac{\sigma^2}{2} \partial_x(x^2 \hat{f}_{hk}(x,t)) \right],
\end{equation}
with 
\[
\mathcal A_{hk\ell r}(x) = - \int_{\mathbb R^2} x  \Phi_\delta(x/x_L,\z)\Psi^1_h(z_1)\Psi^2_k(z_2)\Psi^1_\ell(z_1)\Psi^2_r(z_2) \rho_1(z_1)\rho_2(z_2)dz_1 dz_2. 
\]
Similarly, we can consider the approximated solution of the obtained von Bertalanffy kinetic model 
\[
f(\z,x,t) = f^M(\z,x,t) \approx \sum_{h,k,\ell = 0}^M \hat{f}_{hk\ell}(x,t) \Psi^1_h(z_1)\Psi^2_k(z_2)\Psi^3_\ell(z_3)
\]
that is determined by the following set of equations
\begin{equation}\label{eq:sG_vB}
    \partial_t \hat f_{hk\ell} = \partial_x \left[ \sum_{p,r,s = 0}^{M} \mathcal A_{hk\ell prs} \hat{f}_{p r s}(x,t) + \dfrac{\sigma^2}{2} \partial_x(x^2 \hat{f}_{hk\ell}(x,t)) \right],
\end{equation}
and 
\[
\mathcal A_{hk\ell prs}(x) = - \int_{\mathbb R^3} x  \Phi_\delta(x/x_L,\z)\Psi^1_h(z_1)\Psi^2_k(z_2)\Psi^3_\ell(z_3)\Psi^1_p(z_1)\Psi^2_r(z_2)\Psi^3_s(z_3) \rho_1(z_1)\rho_2(z_2)\rho_3(z_3)dz_1 dz_2dz_3. 
\]

\subsection{Free growth case: convergence and agreement with available data} \label{subsec::conv}

In the following, we show the convergence of the sG scheme for the Fokker-Planck equation \eqref{eq:FP}. We consider deterministic initial conditions coherent with growth curves of Figure \ref{fig:curve_sperimentali} after $100$ days from the  tumour onset. These observations are distributed as a Gamma density
\[
f_{0}(x) = \dfrac{p_1^{p_2} x^{p_1-1} e^{-p_2x}}{\Gamma(p_1)},
\]
with $(p_1,p_2)=(0.3,2.8)$ for the Gompertz case and with  $(p_1,p_2)=(0.37,2.2)$ for the von Bertalanffy case.

We introduce then a uniform discretisation of the domain  $[0,\,2]\subset \mathbb R_+$ obtained with $N=201$ gridpoints, $\Delta x=10^{-2}$ and a time discretisation of the interval $[0,\,T]$ obtained with $\Delta t = \Delta x/C$ with $C = 10^2$ and $T=10$. A central difference scheme is then considered for the numerical solution of the systems of equations \eqref{eq:sG_Gompertz}-\eqref{eq:sG_vB}.

\begin{figure}
	\centering
	\includegraphics[scale = 0.56]{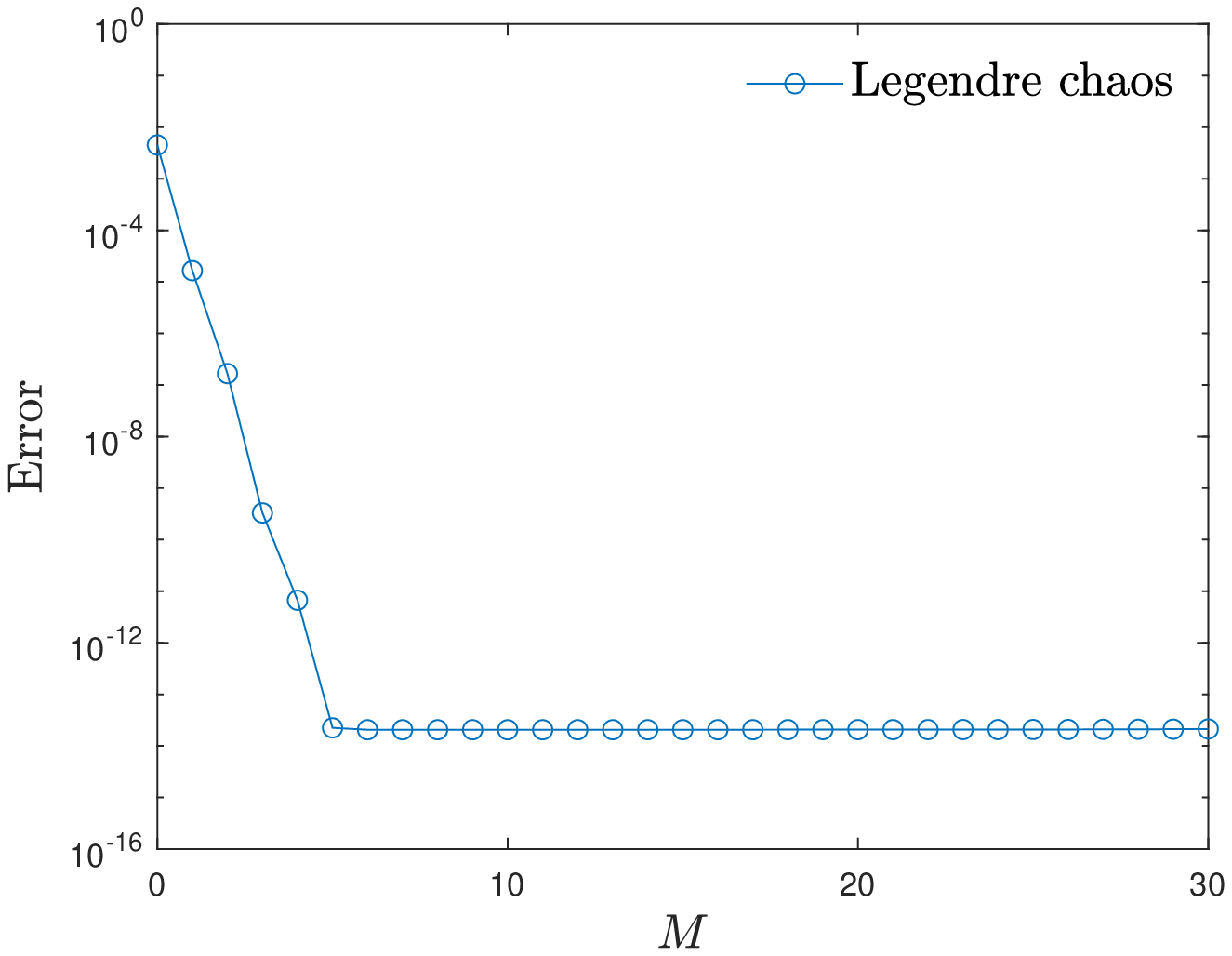}
	\includegraphics[scale = 0.56]{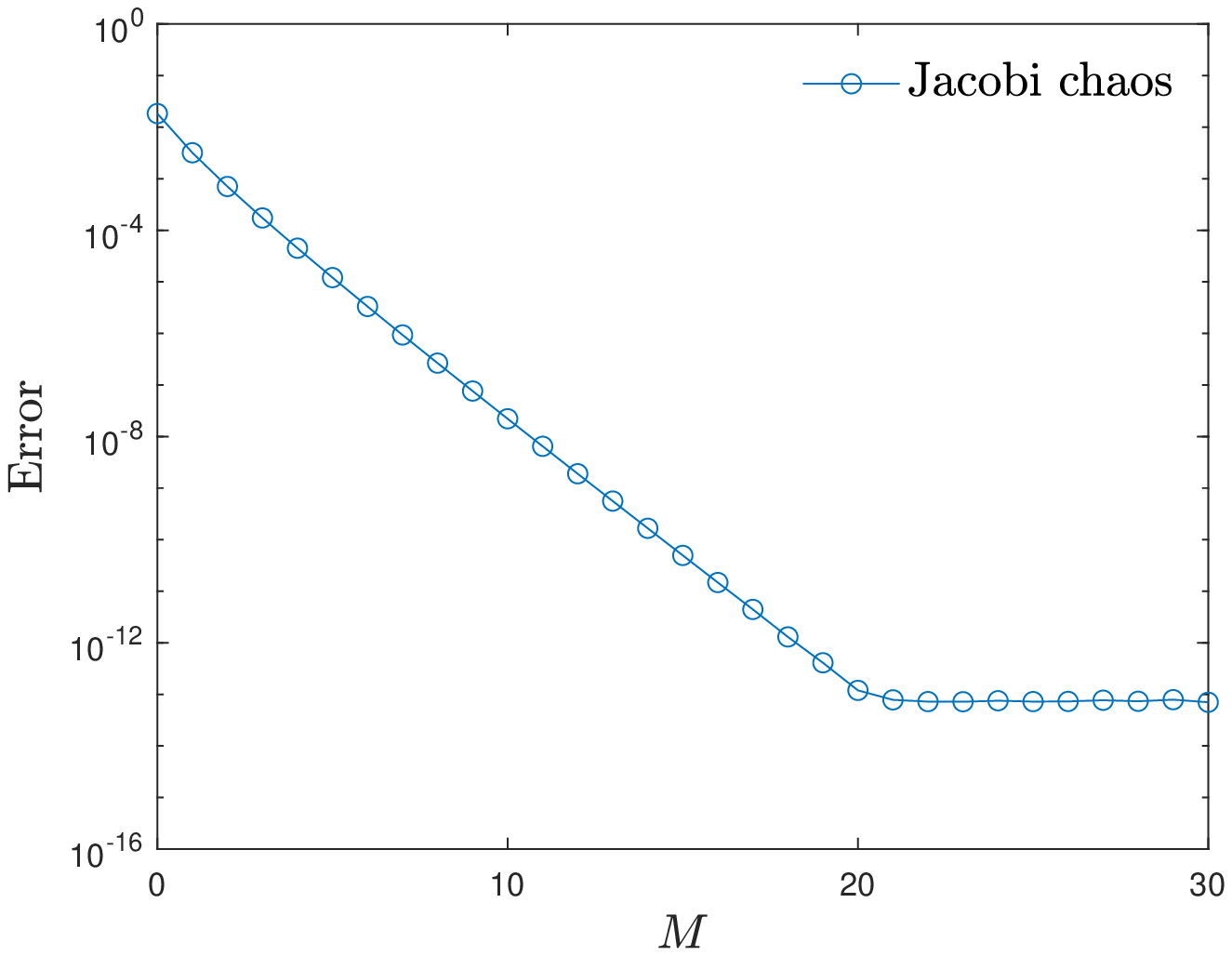} \\
    \includegraphics[scale = 0.56]{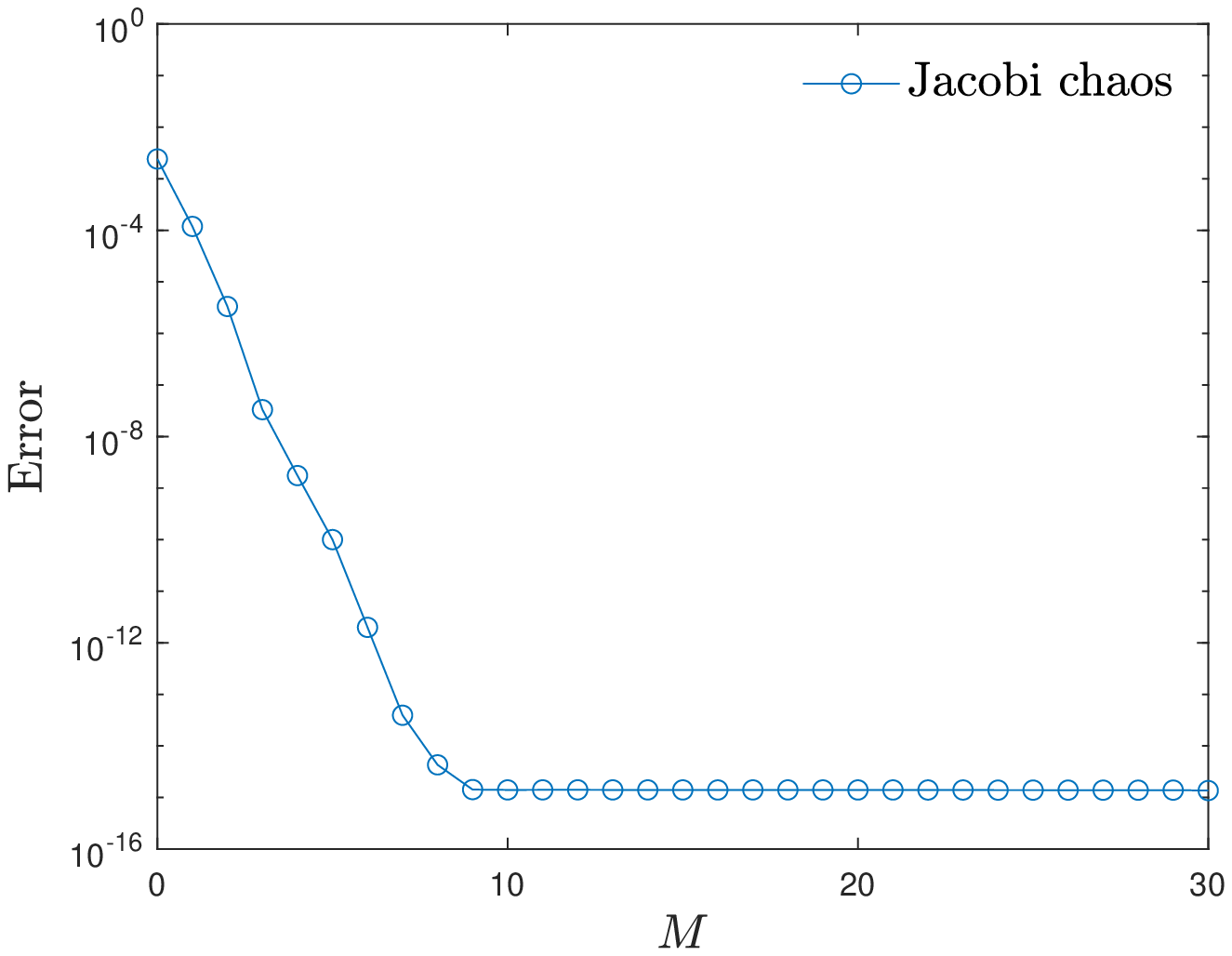}
	\includegraphics[scale = 0.56]{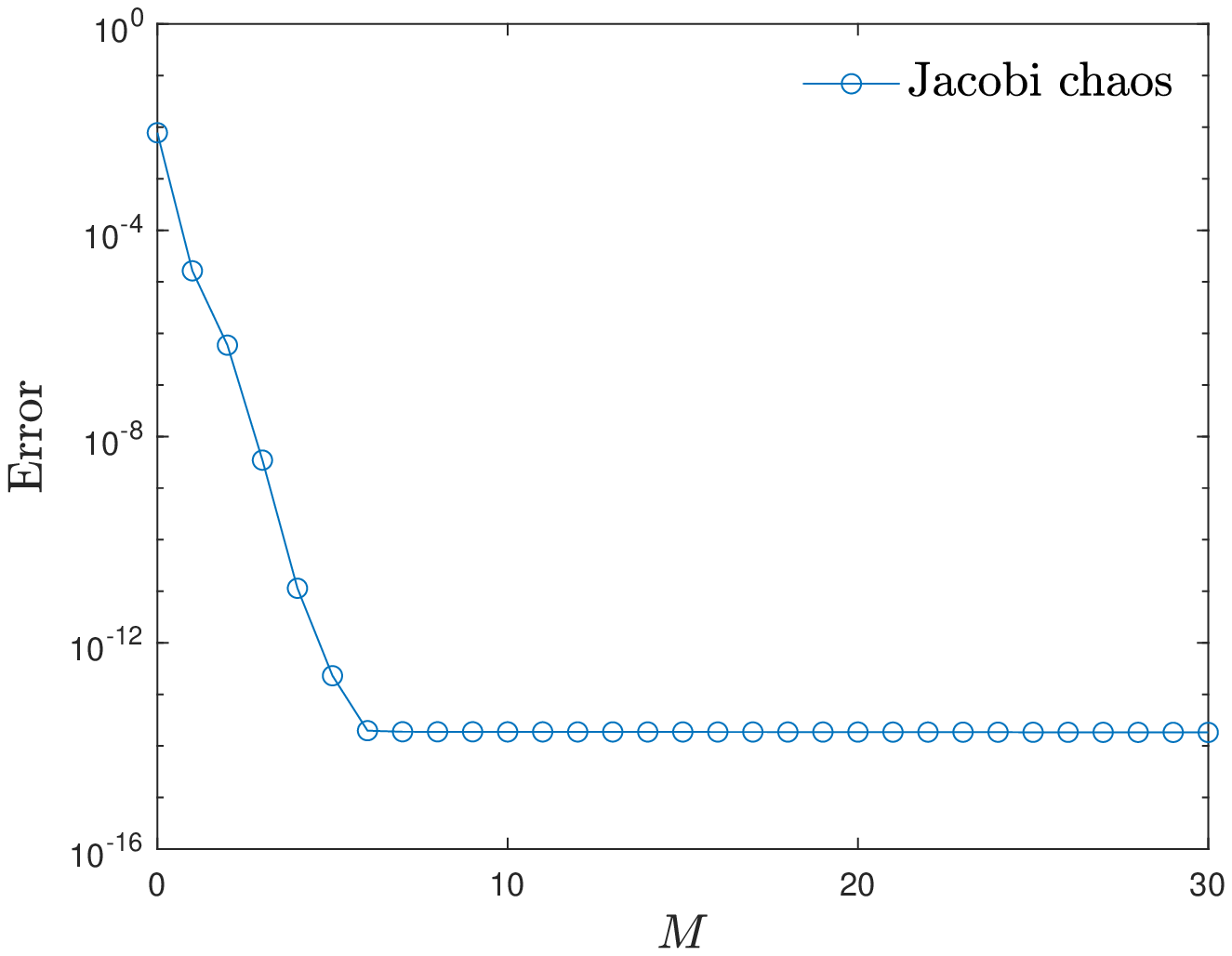}
	\caption{Convergence of the $L^2$ error with respect to a reference solutions obtained with $\bar M = 50$, at fixed time $T=10$. The top row corresponds to the Gompertz case whereas the bottom row to the von Bertalanffy case. Top-Left:  $z_1=\alpha \sim \mathcal U([10^{-3},3 \cdot 10^{-2}])$ and fixed $z_2=x_L\equiv 0.5$. Top-Right: $z_2=x_L\sim \textrm{B}(c_1,c_2)$ and fixed $z_1=\alpha\equiv 0.01$. Bottom-Left:  $z_1=a \sim \textrm{B}(c_1,c_2)$ and fixed $z_2=q\equiv 0.01 $, $z_3=x_L\equiv 0.5$. Bottom-Right: $z_2=q\sim \textrm{B}(c_1,c_2)$ and fixed $z_1=a\equiv 0.8$, $z_3=x_L\equiv 0.5$. The values $(c_1,c_2)$ are reported in Table \ref{table:tabella}.
	}
	\label{fig:test_conv}
\end{figure}

As for the uncertain parameters, we refer to Subsection \ref{subsec:growth}, and in particular to Table \ref{table:tabella}, for the choice of the distributions and, consequently, of the polynomial basis. Let us recall that the Uniform distribution and the Beta distribution corresponds to a Legendre polynomial chaos expansion and a Jacobi polynomial chaos expansion, respectively. 

In the following, we numerically check the convergence of the scheme in the space of random parameters in terms of the evolution of mean volumes. Hence, we consider a reference approximation of the first moment 
\[
m^{\bar M}(\z,t)=\int_{\mathbb R_+} x f^{\bar M}(\z,x,t)dx
\]
obtained with $\bar{M}=50$ at fixed time $T=10$. Then, we compute the $L^2$ error at time $t>0$ defined as 
\[
\|m^{\bar M}(\z,t) - m^M(\z,t) \|_{L^2(\Omega)}
\]
where $m^{M}(\z,T)$ is the first moment obtained with a sG expansion up to order $M<\bar M$, with $M=0,\dots,30$.
\begin{figure}
	\centering
	\includegraphics[scale = 0.56]{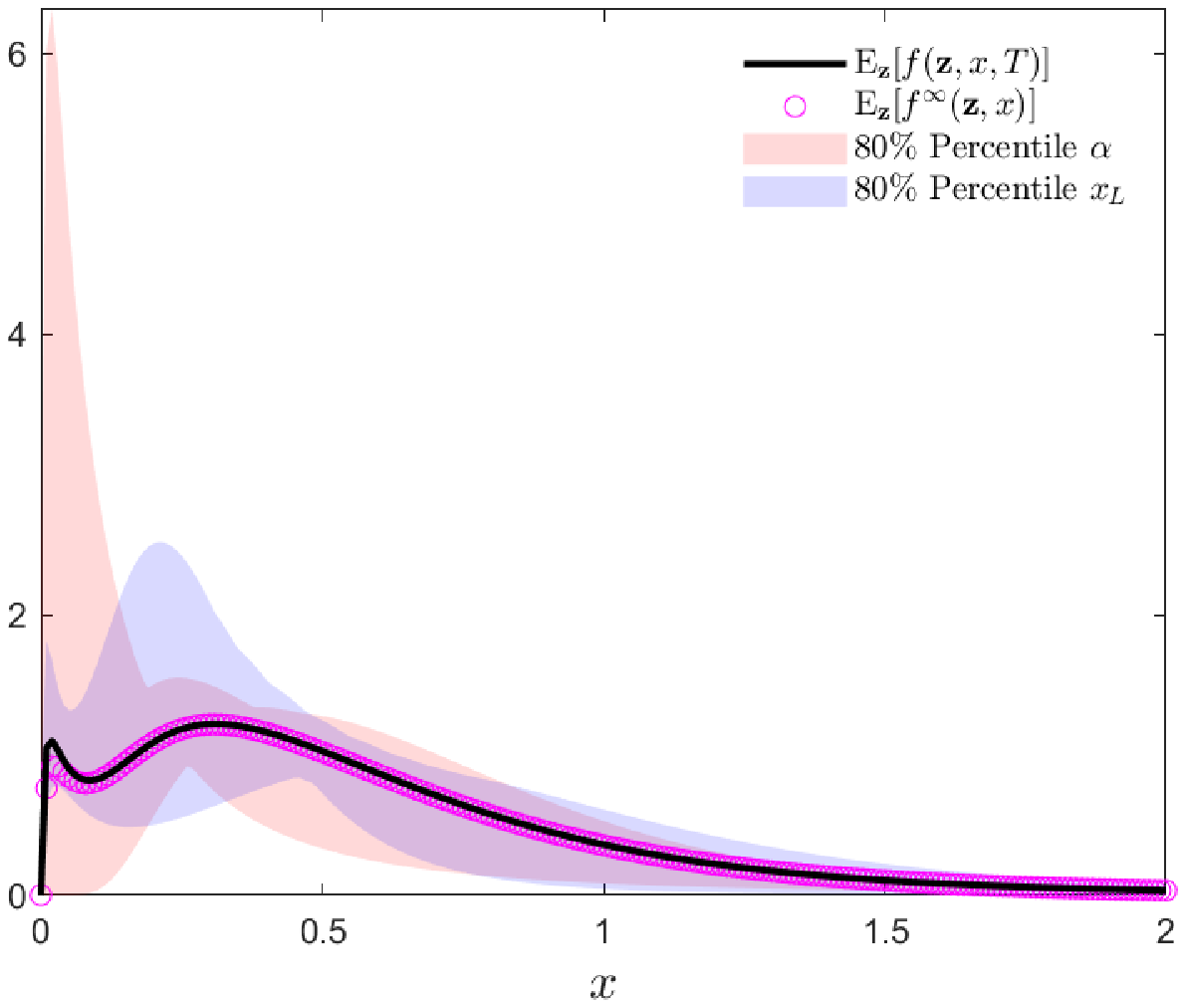}
	\includegraphics[scale = 0.56]{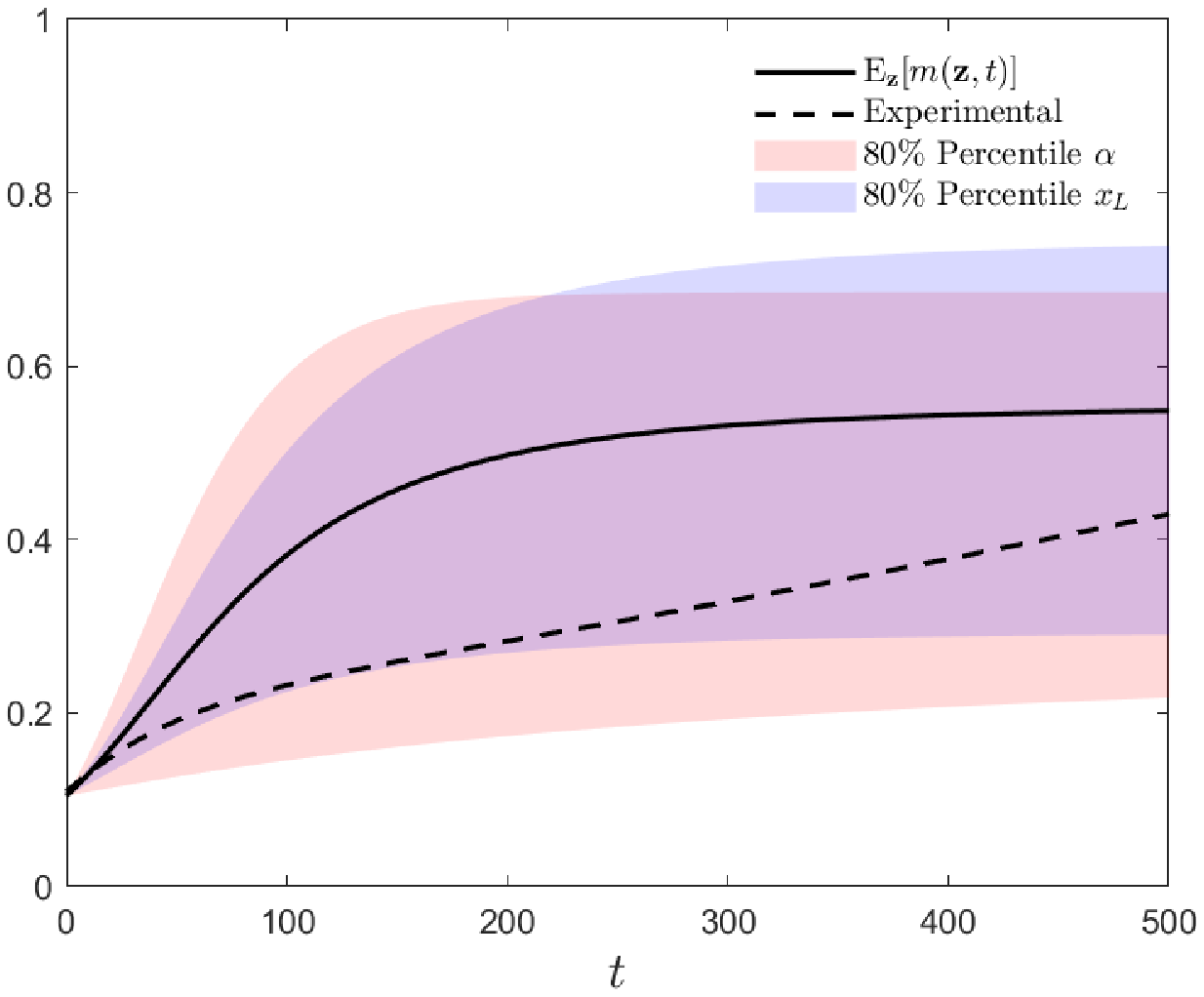} \\
	\includegraphics[scale = 0.56]{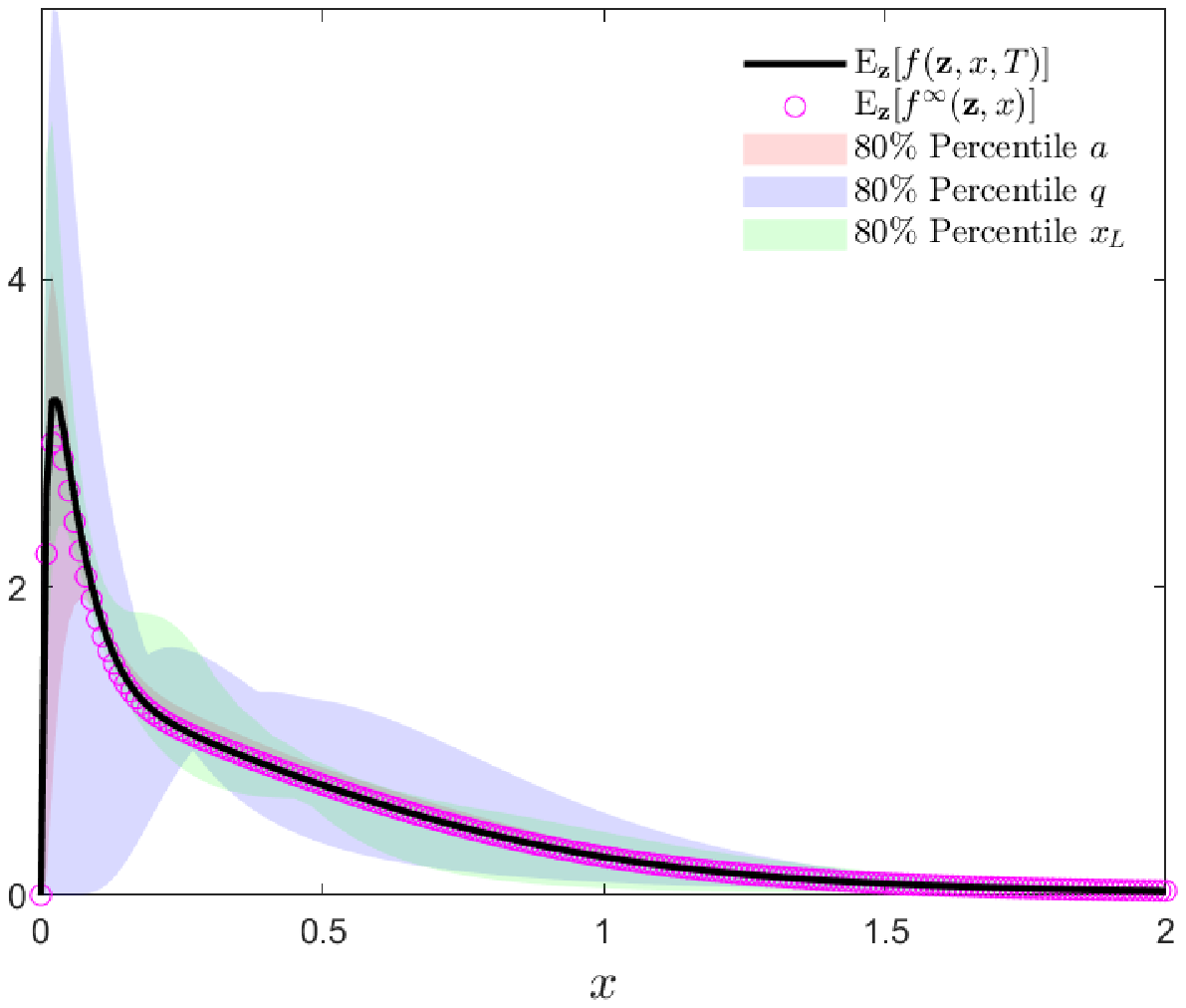}
	\includegraphics[scale = 0.56]{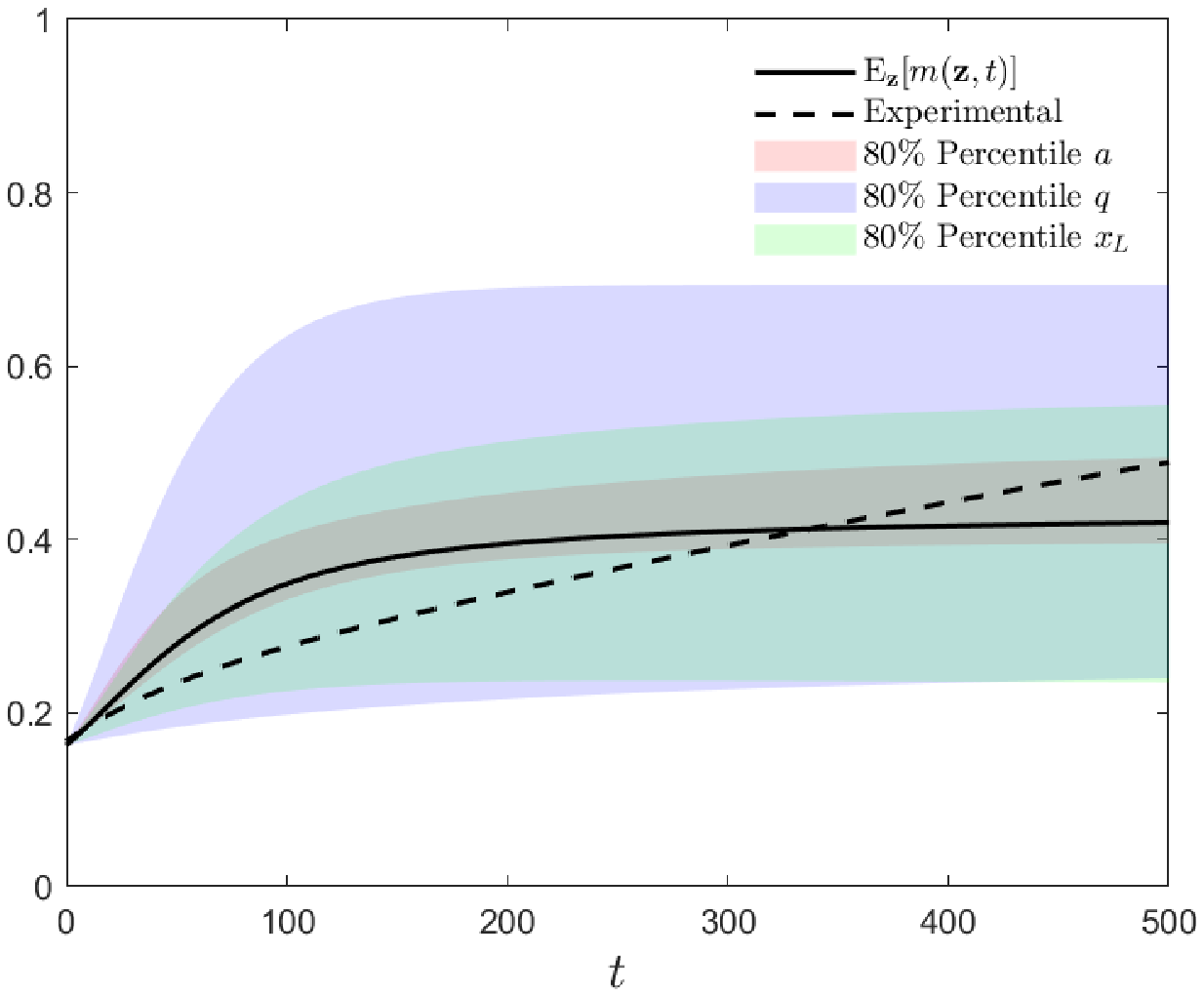}
	\caption{Top: large time distribution (left) and evolution of the mean volume (right) for the Gompertz kinetic model with 2D uncertainties. The solid line is the numerical solution of \eqref{eq:FP} at the final time $T=500$, the markers refers to the expectation of the analytic solution \eqref{eq:equilibria}. Bottom: large time distribution (left) and evolution of the mean volume (right) for the von Bertalanffy kinetic model with 3D uncertainties. In all the cases, we choose $\Delta x = 10^{-2}$, $\Delta t = \Delta x / C$ with $C=10^2$ and $M=3$.}
	\label{fig:test_nc}
\end{figure}
In Figure \ref{fig:test_conv} we may observe the rapid decay of the numerical error in the random space in both models that we have considered.  We observe that we reach essentially the machine precision with a relatively small number of projections.

Once we have checked for the spectral convergence of the method, we can investigate the behaviour of our models with respect to the experimental data. In particular, we will look at the QoI introduced in Section \ref{sec:QoI}. To this aim, we use the introduced numerical setting with $M = 3$ in all the subsequent numerical tests. In the top row of Figure \ref{fig:test_nc} we show the emerging equilibrium distribution from \eqref{eq:sG_Gompertz} with the discussed 2D uncertainty and in the right plot the evolution of the mean volume of the tumours. In the bottom row of Figure  \ref{fig:test_nc} we concentrate on the model \eqref{eq:sG_vB} with 3D uncertainty and again the evolution of the mean volume of tumours in the right plot. The shaded colour bands are relative to the variability computed with percentiles with respect to the introduced uncertain quantities.

\begin{figure} 
    \centering
    \begin{subfigure}[b]{0.475\textwidth}
        \centering
        \includegraphics[width=\textwidth]{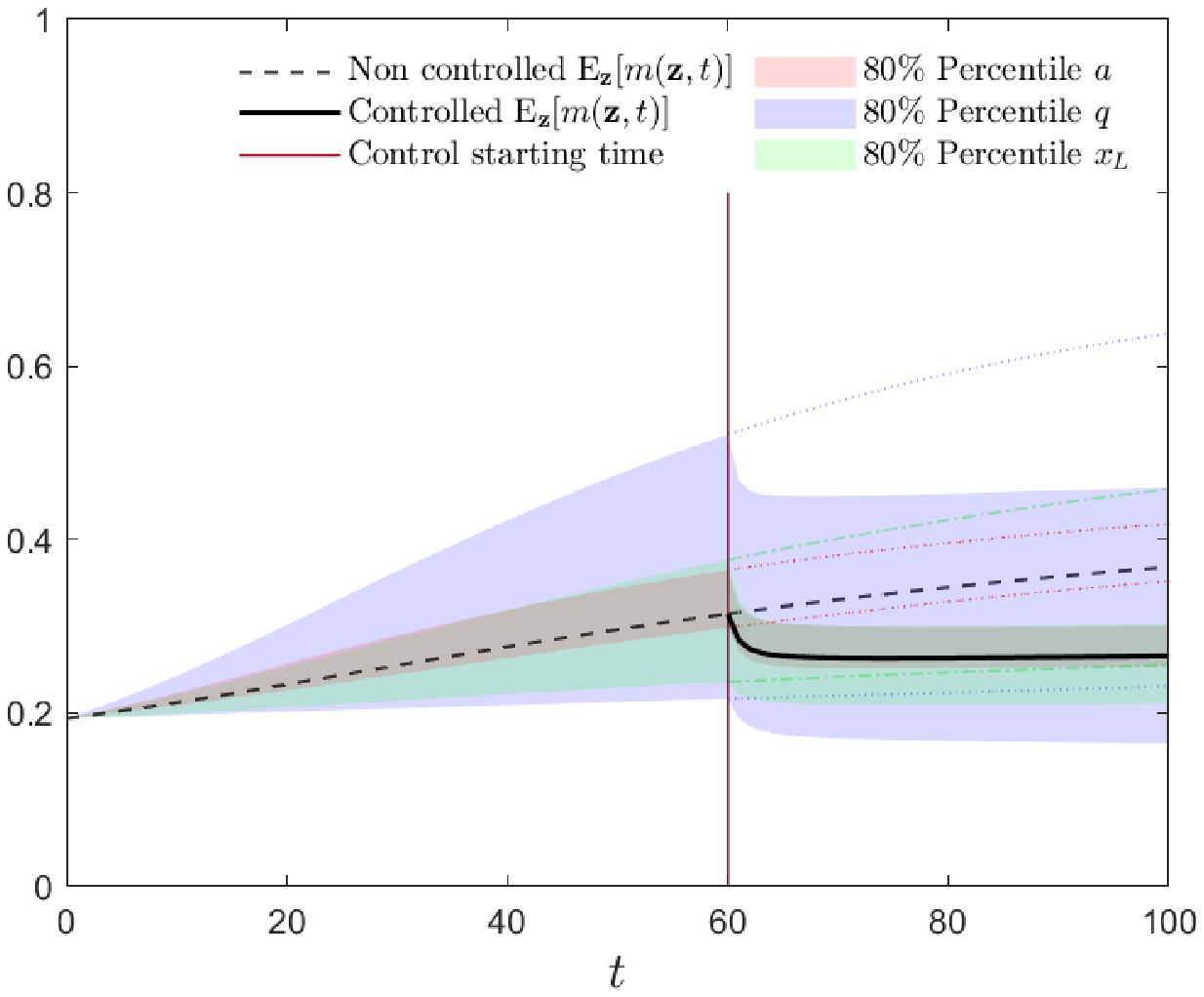}
        \caption[Network2]%
        {{\small $p=1$, $\kappa=1$, $S(x)=1$}}    
    \end{subfigure}
    \hfill
    \begin{subfigure}[b]{0.475\textwidth}  
        \centering 
        \includegraphics[width=\textwidth]{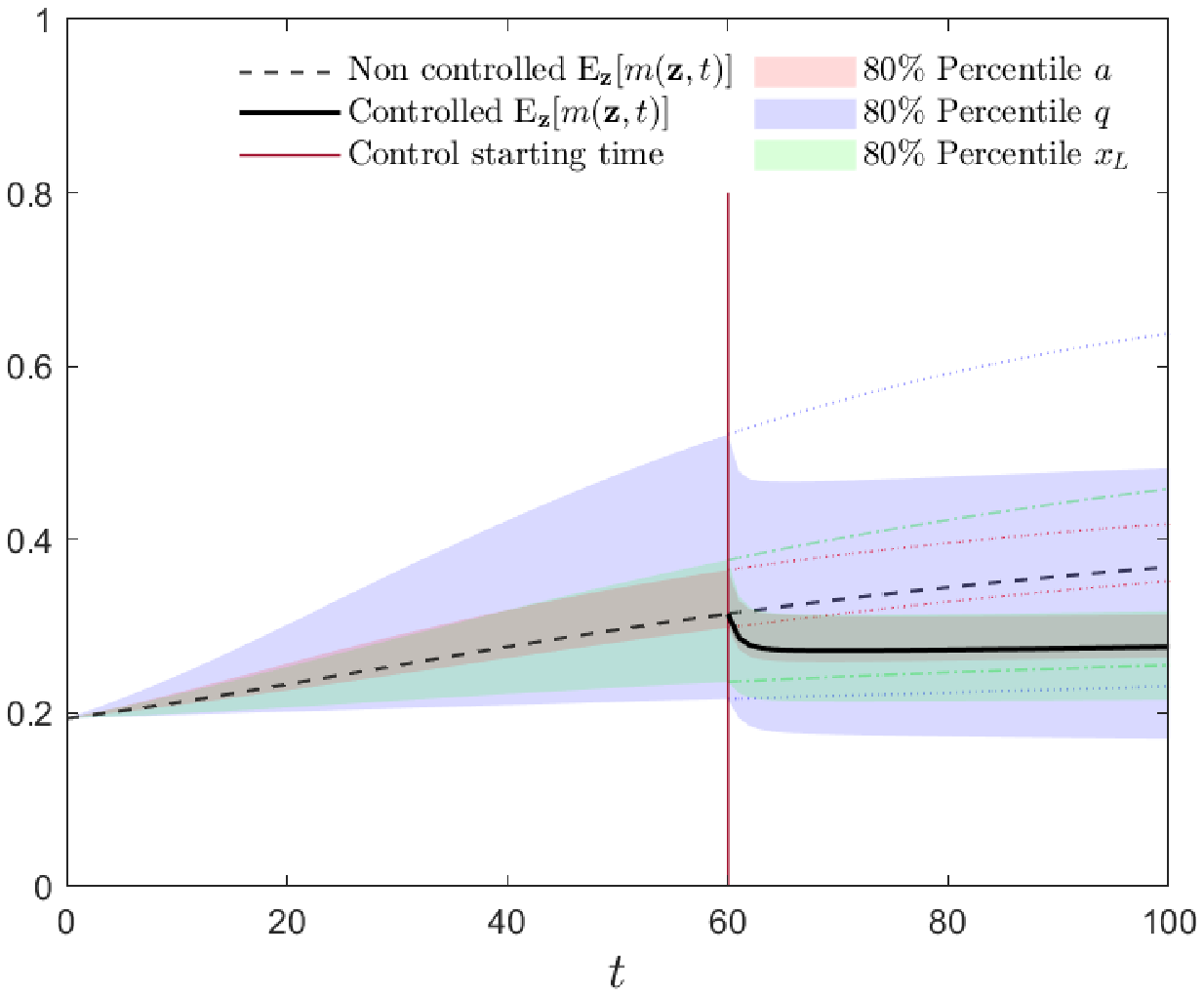}
        \caption[]%
        {{\small $p=1$, $\kappa=1$, $S(x)=\sqrt{x}$}}    
    \end{subfigure}
    %
    \begin{subfigure}[b]{0.475\textwidth}   
        \centering 
        \includegraphics[width=\textwidth]{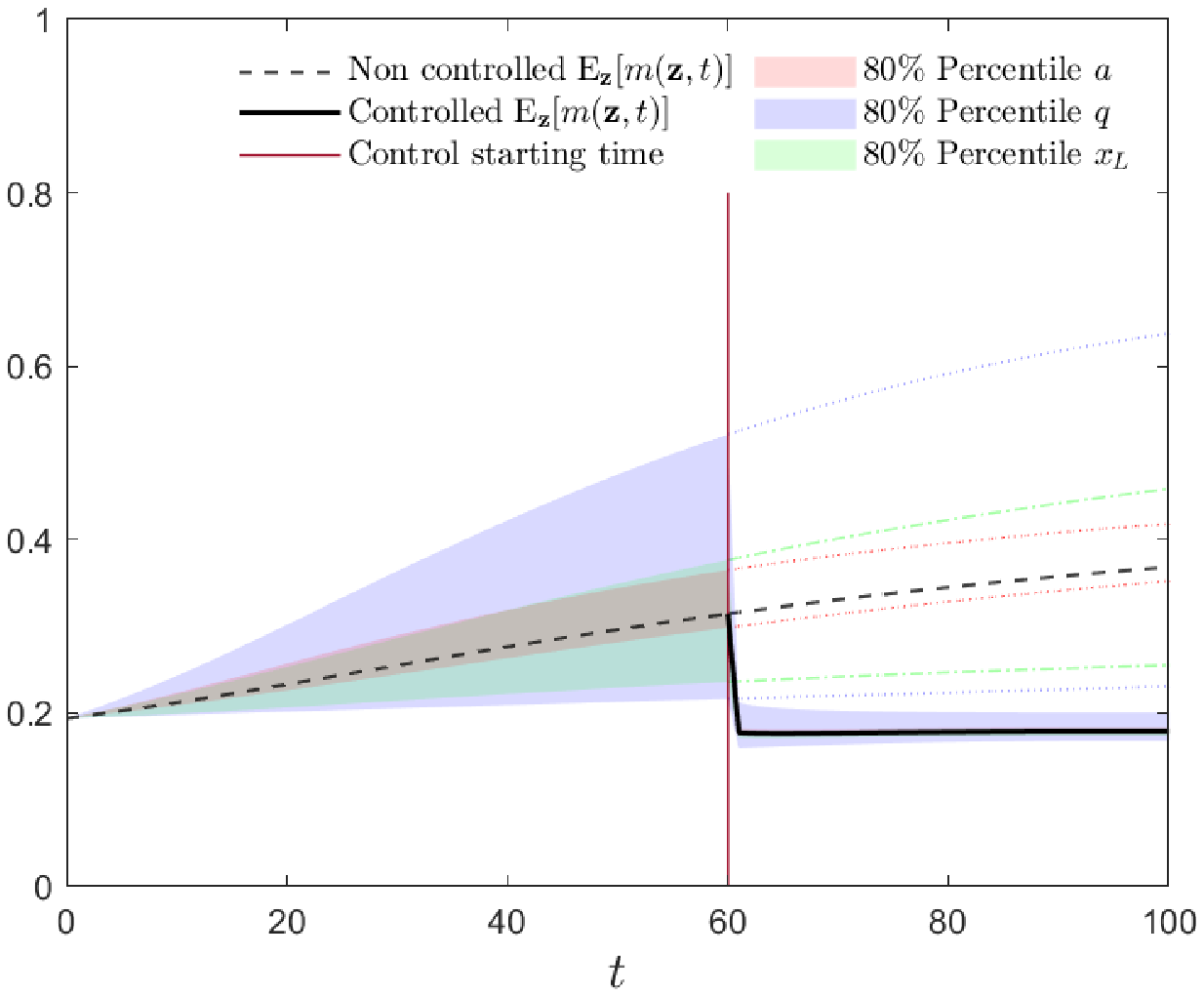}
        \caption[]%
        {{\small $p=1$, $\kappa=0.1$, $S(x)=1$}}    
    \end{subfigure}
    \hfill
    \begin{subfigure}[b]{0.475\textwidth}   
        \centering 
        \includegraphics[width=\textwidth]{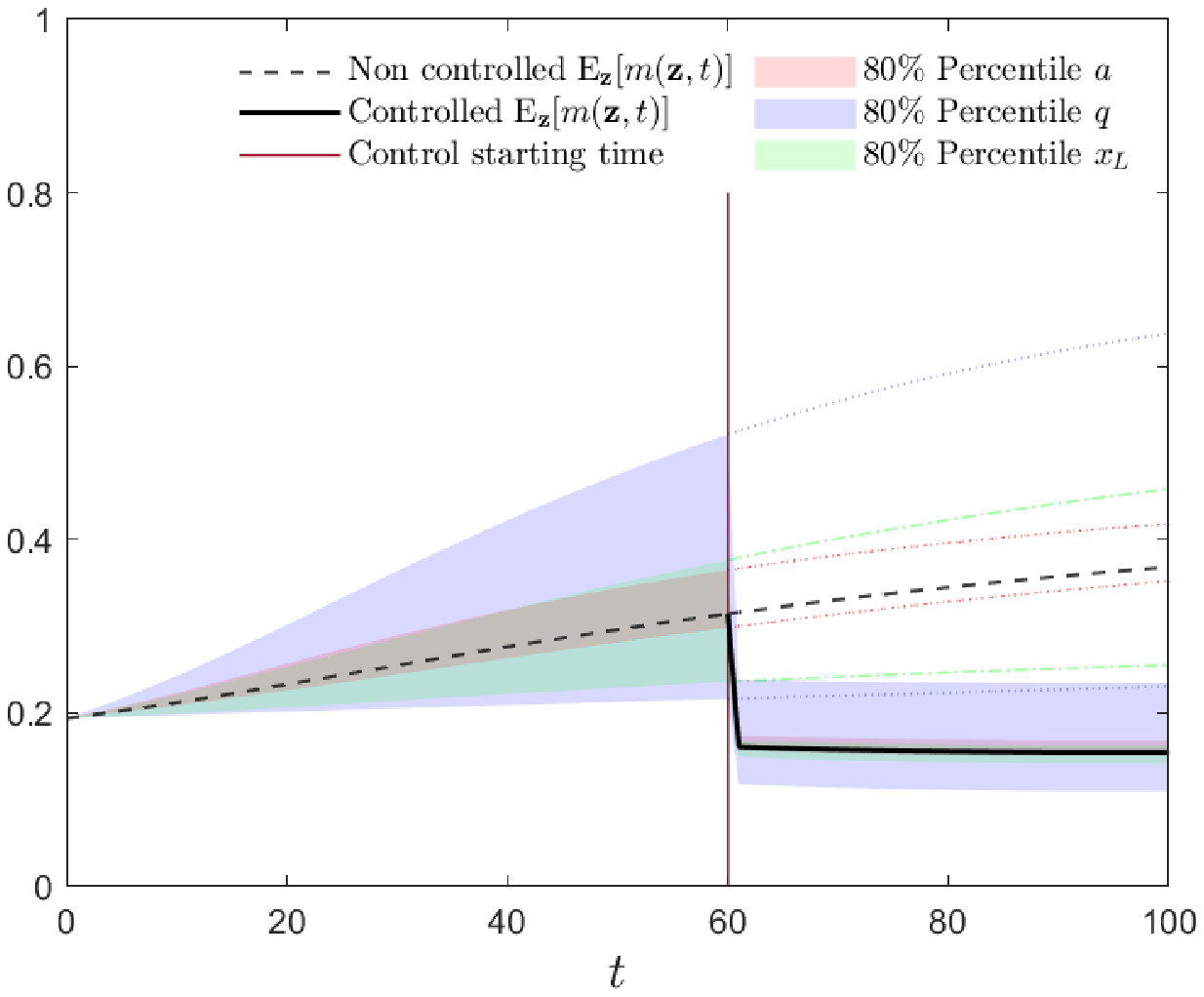}
        \caption[]%
        {{\small $p=1$, $\kappa=0.1$, $S(x)=\sqrt{x}$}}    
    \end{subfigure}
    %
    \begin{subfigure}[b]{0.475\textwidth}   
        \centering 
        \includegraphics[width=\textwidth]{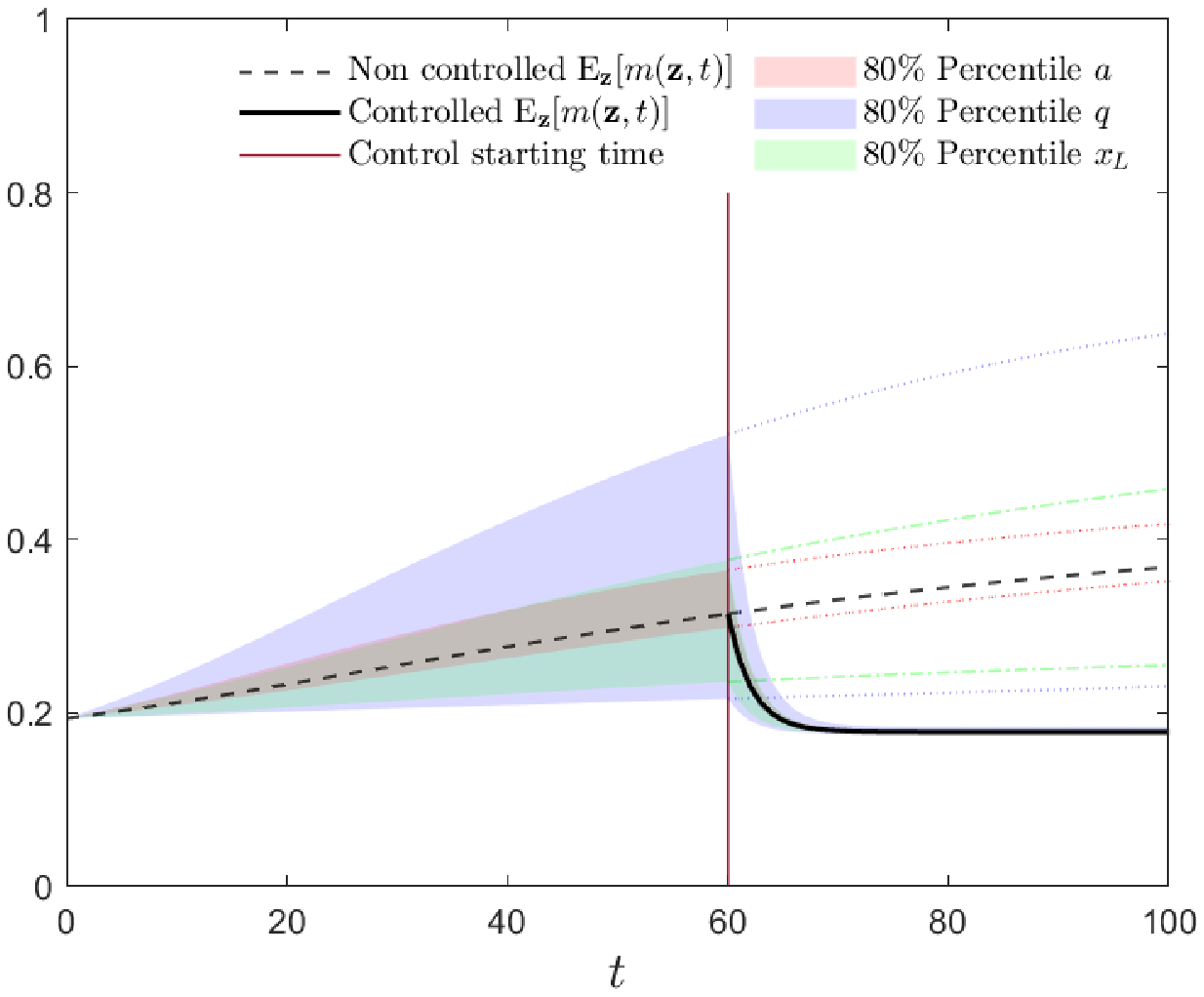}
        \caption[]%
        {{\small $p=2$, $\kappa=1$, $S(x)=1$}}    
    \end{subfigure}
    \hfill
    \begin{subfigure}[b]{0.475\textwidth}   
        \centering 
        \includegraphics[width=\textwidth]{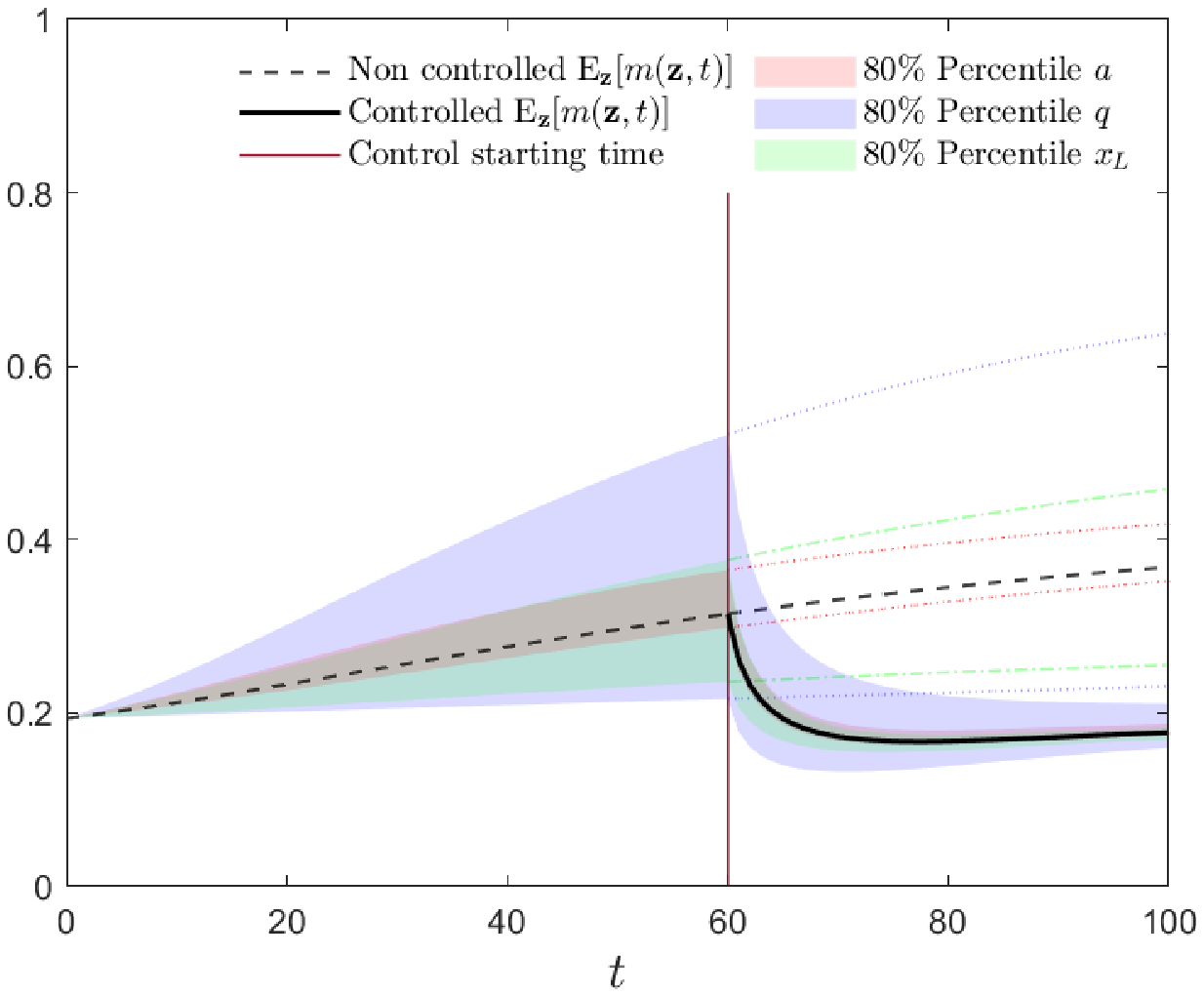}
        \caption[]%
        {{\small $p=2$, $\kappa=1$, $S(x)=\sqrt{x}$}}    
    \end{subfigure}

    \caption{Evolution of $m(\z,t)$ in the uncontrolled scenario for $t\leq60$ and in a controlled scenario for $t>60$, with $p=1,2$, $S(x)=1,\sqrt{x}$ and $\kappa=0.1,1$. The dashed lines represent the trend of the tumour if the control is not in action. In the case $p=2$, we adopted a sG scheme with $\Delta x=10^{-2}$, $\Delta t=\Delta x/C$ with $C=10^2$, $M=3$ for the numerical solution of \eqref{eq:sG_vB} with the introduced clinical uncertainties. In the case $p=1$, we adopted a stochastic collocation DSMC with $N=10^5$, $\Delta t=0.05$, $\epsilon=2\Delta t$ and $M=3$. We considered the experimental target volume size {$x_d=0.18 \times 10^5 mm^3$ and results are scaled by a factor $10^5$}.} 
    \label{fig:VB_p2_m}
\end{figure}

\subsection{Effects of the control and damping of uncertainties}

We consider now the kinetic models \eqref{eq:FPcontr} in presence of control strategies to test the effectiveness of the introduced control  in reducing the tails of the distributions and damping the uncertainties of the system. For this reason, we consider here only the von Bertalanffy case that induces a power-law-type equilibrium distribution as discussed in Section \ref{subsec:equilibria}. From experimental measurements we observed an average value of the target volume $x_d = 0.18\times 10^5 \textrm{mm}^3$ after the implemented therapeutic protocols, for this reason we have fixed this value in each  experiment of this section. {The obtained value of the target volume will be scaled by a factor $10^5$ through the section. To activate the control we compute the mean tumours' size from experimental data. Hence, $u$ starts acting when $\mathbb E_{\z}[m(\z, t)]$ exceeds this threshold.}

In Figure \ref{fig:VB_p2_m} we present the evolution of the expected values of the first order moment $m(\z,t)$ in a constrained setting obtained from \eqref{eq:FPcontr}. In particular, we plot the uncontrolled evolutions up to the time $t=60$ and then we activate the control. As in the uncontrolled scenario, we consider a uniform discretisation of $[0,\,2]$ obtained with $N=201$ gridpoints and a time step $\Delta t=\Delta x/C$, with $C=10^2$, for the time interval $[0,\,T]$, with $T=100$ final time. We notice that the control succeeds in reducing both the expected values of the first moment and the uncertainty, with smaller values of $\kappa$.

\begin{figure}
    \centering
    \begin{subfigure}[b]{0.475\textwidth}
        \centering
        \includegraphics[width=\textwidth]{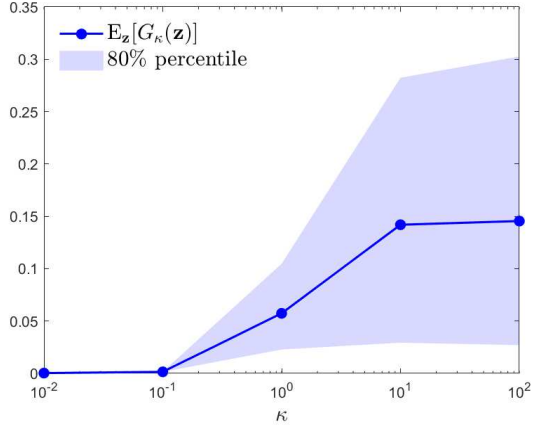}
        \caption[Network2]%
        {{\small $p=1$, $S(x)=1$}}    
    \end{subfigure}
    \hfill
    \begin{subfigure}[b]{0.475\textwidth}  
        \centering 
        \includegraphics[width=\textwidth]{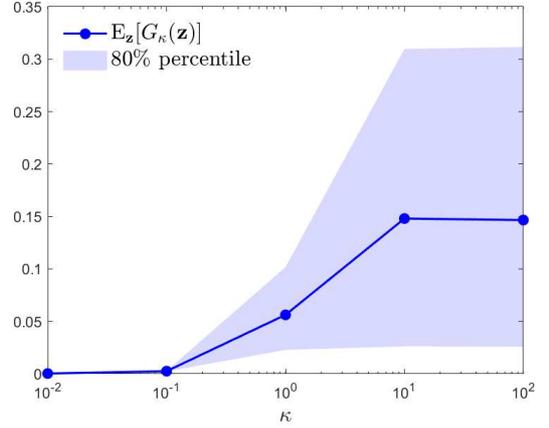}
        \caption[]%
        {{\small $p=1$, $S(x)=\sqrt{x}$}}    
    \end{subfigure}
    \vskip\baselineskip
    \begin{subfigure}[b]{0.475\textwidth}   
        \centering 
        \includegraphics[width=\textwidth]{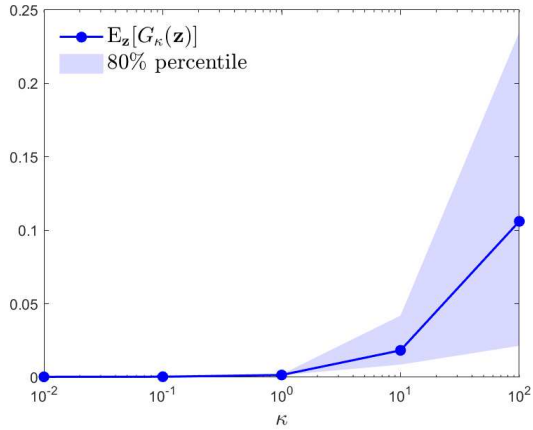}
        \caption[]%
        {{\small $p=2$, $S(x)=1$}}    
    \end{subfigure}
    \hfill
    \begin{subfigure}[b]{0.475\textwidth}   
        \centering 
        \includegraphics[width=\textwidth]{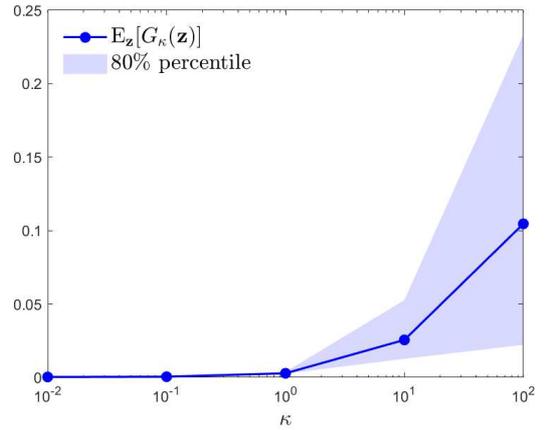}
        \caption[]%
        {{\small $p=2$, $S(x)=\sqrt{x}$}}    
    \end{subfigure}
    \caption[ ]
    {Expectation of the quantity $G_{\kappa}(\z)$ defined in \eqref{eq:G} {and computed with $x_d = 0.18$} versus the penalisation $\kappa$, for $p=1,2$ and $S(x)=1,\sqrt{x}$, considering the von Bertalanffy model with uniform-distributed coefficients. The plots are in semi-logarithmic scale to highlight the uncertainty damping for small values of $\kappa$. In all the cases, a collocation DSMC scheme is adopted, with the choices $N=10^5$, $M=3$, $\Delta t=0.05$ and $\epsilon=2\Delta t$.} 
     \label{fig:G_p1p2}
\end{figure}

To quantify the effectiveness of the adopted control strategy, we define an index that quantifies the variability around the target $x_d$ computed at a given time $T>0$ and defined as follows
\be \label{eq:G}
G_{\kappa}(\z)=\int_{\mathbb{R}_+}(x-x_d)^2 f(\z,x,T)dx,
\ee
where $f(\z,x,T)$ is the kinetic distribution of the controlled model with embedding the penalisation coefficient $\kappa>0$. In Figure \ref{fig:G_p1p2} we show the behaviour of $\mathbb E_\z[G_{\kappa}(\z)]$ together with confidence bands and computed for several penalisation coefficients. We considered both the cases $p=1,2$ and selective functions $S(x) = 1,\sqrt{x}$. 

For the case $p = 1$ we adopt a stochastic collocation approach for the kinetic model \eqref{eq:kinetic_ther_strong} that is solved through a DSMC scheme \cite{TZ0}. We choose $N=10^5$, $M=3$ and $\Delta t=0.05$ and $\epsilon=2\Delta t$. We notice that, in all the considered cases, $\mathbb{E}_{\z}[G_{\kappa}(\z)]$ decreases with smaller values of $\kappa$ and the uncertainty is {dampened}.

\begin{figure}
    \centering
    \begin{subfigure}[b]{0.475\textwidth}
        \centering
        \includegraphics[width=\textwidth]{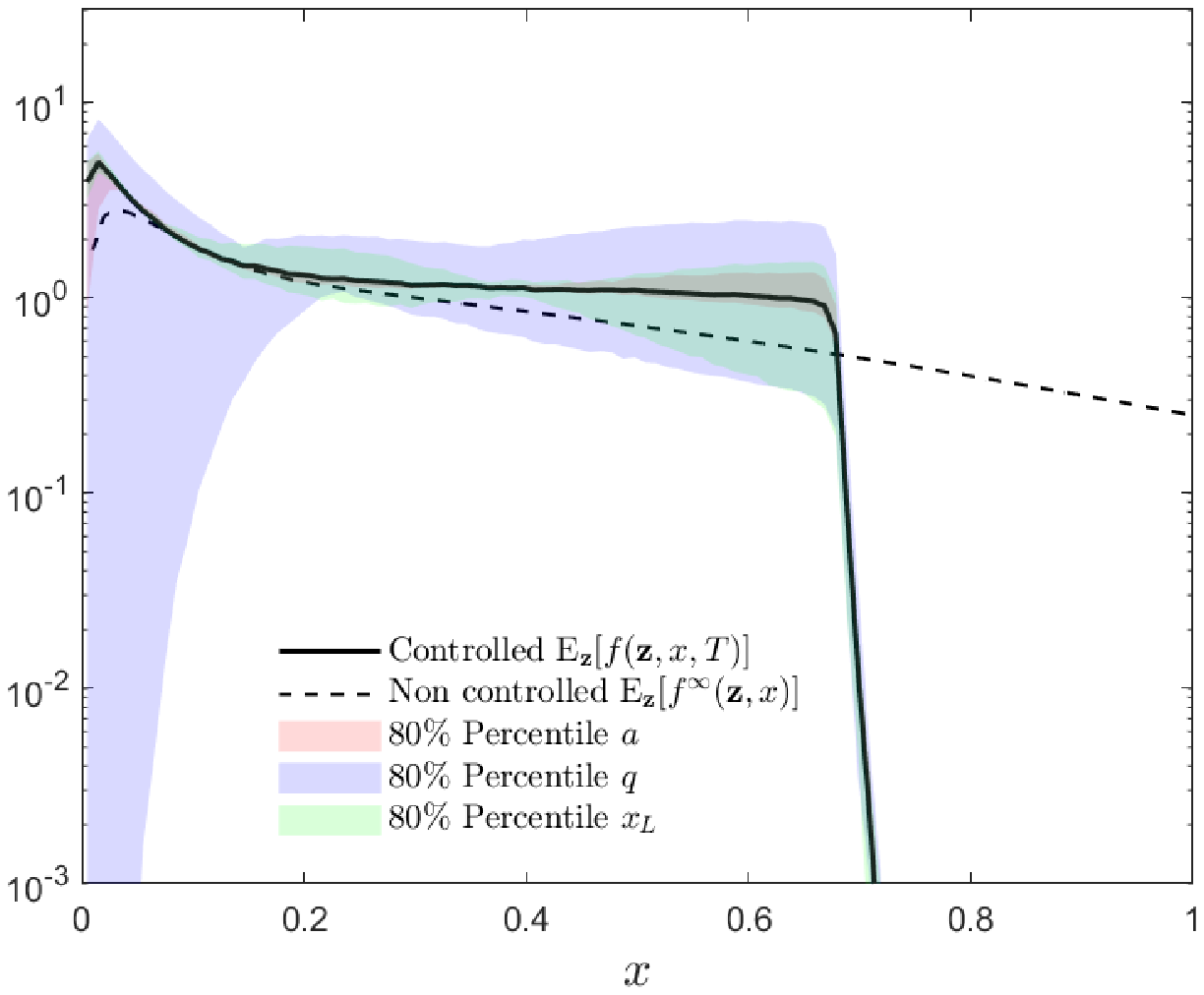}
        \caption[Network2]%
        {{\small $p=1$, $\kappa=1$, $S(x)=1$}}    
    \end{subfigure}
    \hfill
    \begin{subfigure}[b]{0.475\textwidth}  
        \centering 
        \includegraphics[width=\textwidth]{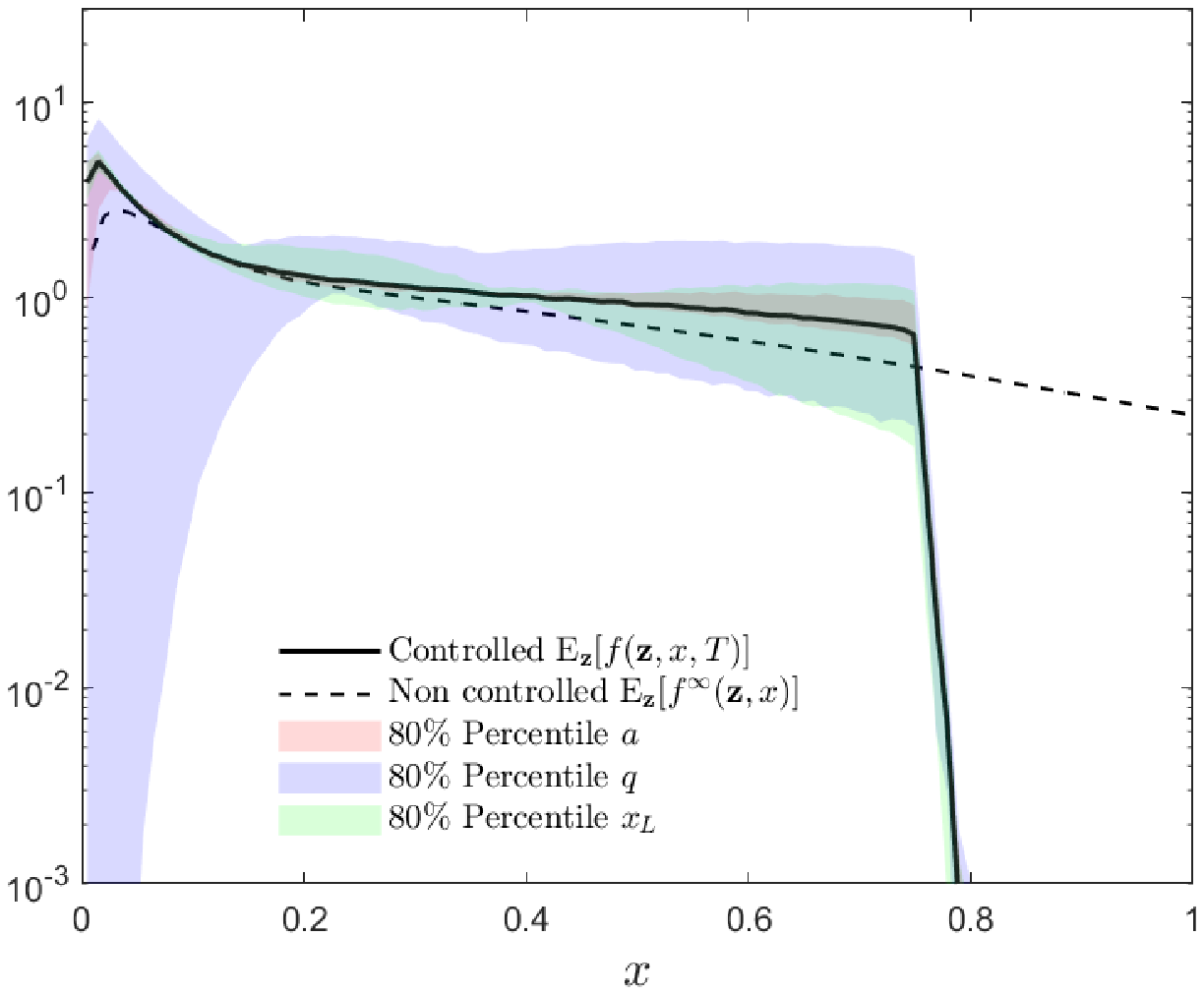}
        \caption[]%
        {{\small $p=1$, $\kappa=1$, $S(x)=\sqrt{x}$}}    
    \end{subfigure}
    \vskip\baselineskip
    \begin{subfigure}[b]{0.475\textwidth}   
        \centering 
        \includegraphics[width=\textwidth]{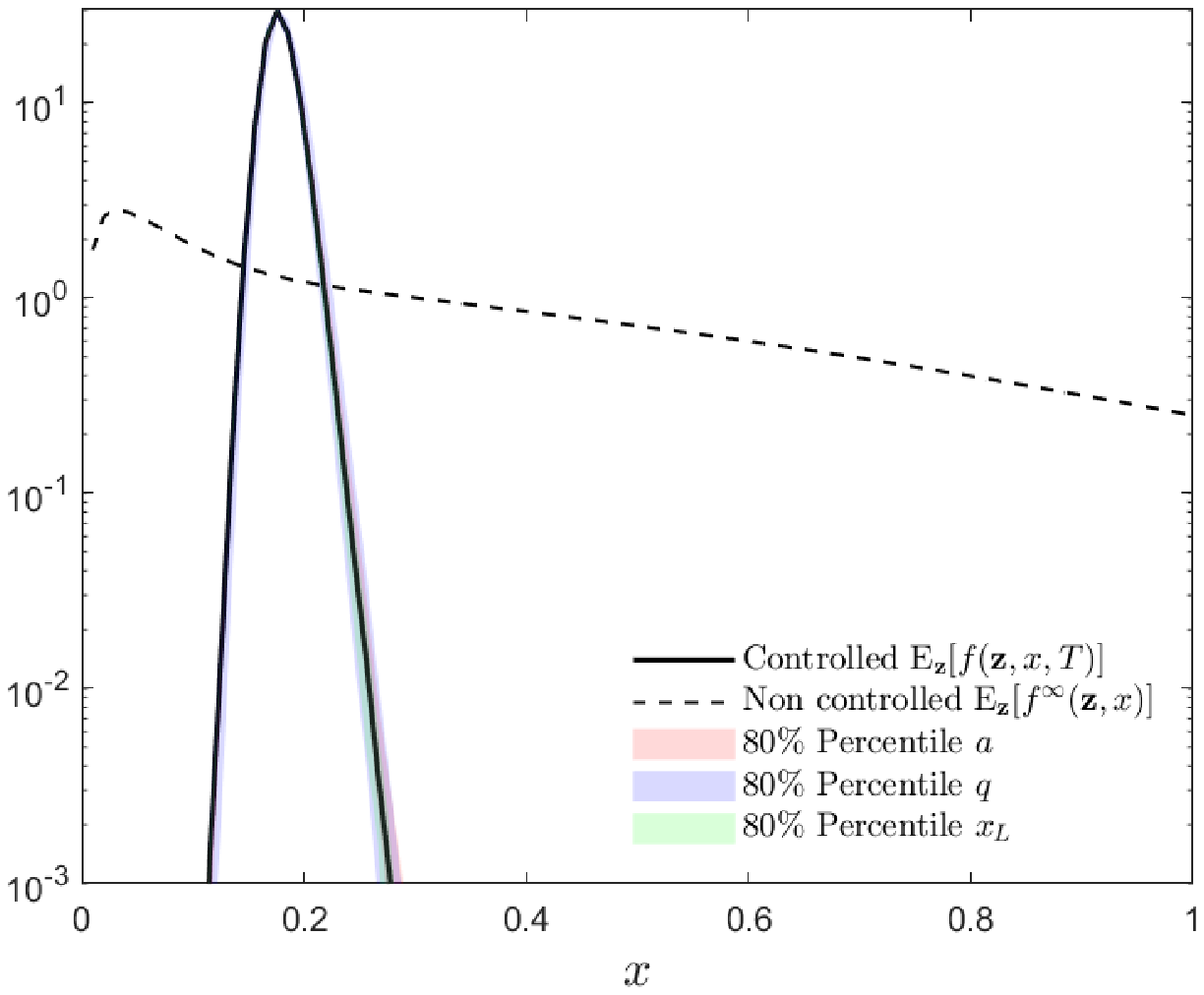}
        \caption[]%
        {{\small $p=2$, $\kappa=1$, $S(x)=1$}}    
    \end{subfigure}
    \hfill
    \begin{subfigure}[b]{0.475\textwidth}   
        \centering 
        \includegraphics[width=\textwidth]{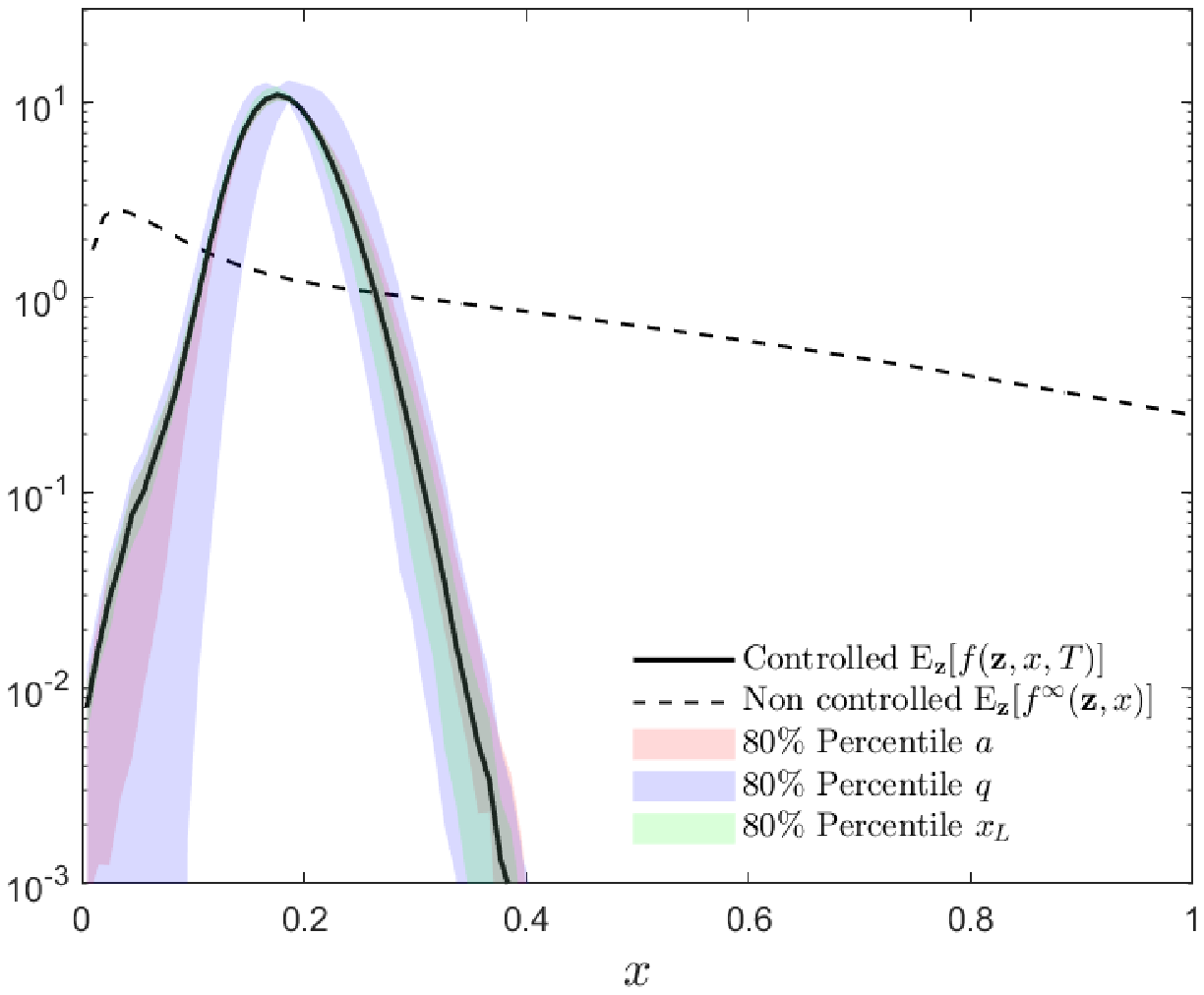}
        \caption[]%
        {{\small $p=2$, $\kappa=1$, $S(x)=\sqrt{x}$}}    
    \end{subfigure}
    \caption[ ]
    {Details of the expected value of the controlled and uncontrolled distributions in the von Bertalanffy growth scenario, for fixed penalisation $\kappa=1$, $p=1,2$ and $S(x)=1,\sqrt{x}$. The solid line is the numerical controlled solution at the final time $T=100$, the dashed line is the uncontrolled analytic solution at the equilibrium. The plots are in semi-logarithmic scale to highlight the suppression of the fat tails. {We considered the experimental target volume size $x_d=0.18 \times 10^5 mm^3$ and results are scaled by a factor $10^5$}} 
     \label{fig:VB_p2_f}
\end{figure}

Now we look directly to the effectiveness of the control strategies in reducing both the tails and the uncertainty of the distributions. In Figure \ref{fig:VB_p2_f} we show the expected distributions in the controlled case for large times, obtained with the introduced sG scheme for the kinetic von Bertalanffy model. We may observe how the introduction of selectivity is capable to change the behaviour of the tail of the distribution as discussed in Section \ref{subsec:FPA}. 

\section*{Conclusions}
In the present paper, we concentrated on the influence of uncertain quantities on kinetic models for tumour growths. Under suitable assumptions, we derived surrogate Fokker-Planck models from which we obtain analytical insight on the large time behaviour of the system. Hence, we proposed suitable selective control strategies mimicking the effects of therapies in steering the volume of tumours towards a target value $x_d$. Through explicit computations, we showed that the solution of the controlled model is close to the target volume and the distance of the first order moment from $x_d$ depends on the penalisation of the control. These control protocols are capable to dampen the variability of the tumours' dynamics due to the presence of uncertainties. Since from the pathological point of view fat-tailed distributions are related to a higher probability of finding large tumours with respect to thin-tailed distributions we observed that by implementing suitable selective strategies we can also change the nature of the emerging distribution of tumours' sizes. Thanks to real observations on a cohort of subjects we observed great variability in the choice of parameters of the models that has been considered in the numerical section. Numerical schemes for the uncertainty quantification of kinetic equations are then considered to observe the effects of the control on the solution of the models. 

\section*{Acknowledgments}
This work has been written within the activities of GNFM group of INdAM (National Institute of High Mathematics). The research was partially supported by Dipartimenti di Eccellenza Program (2018–2022) - Department of Mathematics “F. Casorati”, University of Pavia. M.Z. acknowledges partial support of MUR-PRIN2020 Project (No. 2020JLWP23) ”Integrated mathematical approaches to socio-epidemiological dynamics”. M.Z. would like to thank Prof. Giuseppe Toscani (University of Pavia, Italy) for fruitful and inspiring discussions.

\end{document}